\RequirePackage{fix-cm}
\documentclass[smallextended]{svjour3}       % onecolumn (second format)
\smartqed  % flush right qed marks, e.g. at end of proof
\usepackage{graphicx}
\usepackage{amsmath}
\usepackage{caption}
\usepackage{amssymb}
\usepackage{multicol}
\usepackage{multirow}
\usepackage{epstopdf}
\usepackage{soul}
\usepackage{xcolor}
\usepackage{float}
\usepackage{enumitem}
\graphicspath{{figure/}}

\begin{document}

\title{ENO-Based High-Order Data-Bounded and Constrained Positivity-Preserving Interpolation
}

%\author{First Author         \and
%        Second Author %etc.
%}
\author{T. A. J. Ouermi \\
        %\and
        Robert M. Kirby \\
        %\and
        Martin Berzins \\
        }


\institute{T. A. J. Ouermi \at
           University of Utah School of Computing \\
           U of U Scientific Computing and Imaging Institute \\
           Salt Lake City, Utah, USA \\
           \email{touermi@cs.utah.edu}\\
           \and
           Robert M. Kirby \at 
           University of Utah School of Computing \\
           U of U Scientific Computing and Imaging Institute \\
           Salt Lake City, Utah, USA \\
           \email{kiby@cs.utah.edu}\\
           \and
           Martin Berzins \at
           University of Utah School of Computing \\
           U of U Scientific Computing and Imaging Institute \\
           Salt Lake City, Utah, USA \\
           \email{mb@sci.utah.edu}\\
        }

\date{Received: date / Accepted: date}

\maketitle

\begin{abstract}
  A number of key scientific computing applications that are based upon tensor-product grid constructions, 
  such as numerical weather prediction (NWP) and combustion simulations,  
  require {\em property-preserving} interpolation. Essentially Non-Oscillatory (ENO) interpolation is a
  classic example of such interpolation schemes. 
  In the aforementioned application areas, property preservation often manifests 
  itself as a requirement for either data boundedness or positivity preservation.  For example, in NWP, one may have to interpolate
  between the grid on which the dynamics is calculated to a grid on which the physics is calculated (and back).
  Interpolating density or other key physical quantities without accounting for property preservation may lead to 
  negative values that are nonphysical and result in inaccurate representations and/or interpretations of the physical data.
  
  Property-preserving interpolation is straightforward when used in the context of low-order numerical simulation methods.  High-order property-preserving
  interpolation is, however, nontrivial, especially in the case where the interpolation points are not equispaced. 
  In this paper, we demonstrate that it is possible to construct high-order interpolation methods that ensure either data boundedness or constrained positivity preservation.  
   A novel feature of the algorithm is that the positivity-preserving interpolant is constrained; that is, the amount by which it exceeds the data values may be strictly controlled.
  The algorithm we have developed comes with theoretical estimates that provide sufficient conditions for data boundedness and constrained positivity preservation.
  We demonstrate the application of our algorithm on a collection of 1D and 2D numerical examples, and show that in all cases property preservation 
  is respected.
%subject classification numbers as needed.
  %\PACS{PACS code1 \and PACS code2 \and more}
\end{abstract}

\section{Introduction}
\label{sec:introduction}

A number of key scientific computing applications that are based upon high-order methods over tensor-product grid constructions, 
such as numerical weather prediction (NWP) and combustion simulations, require {\em property-preserving} interpolation.  
In the aforementioned application areas, property preservation often manifests 
itself as a requirement for either data boundedness or positivity preservation. 
The particular application motivating this work is the Navy Environmental Prediction System Using the NUMA Core (NEPTUNE).
NEPTUNE is a next-generation global NWP system being 
developed at the Naval Research Laboratory (NRL) and the Naval Postgraduate School (NPS) \cite{Neptune:Alex}.
NEPTUNE makes use of the Nonhydrostatic Unified Model of the Atmosphere (NUMA) \cite{NeptuneNUMA} three-dimensional spectral element 
dynamical core, but currently uses physics routines that were developed assuming uniform grid spacing on the elements. At least two options are
available for combining these two NWP building blocks:  either (1) evaluate the physics routines at the (nonuniformly-spaced) quadrature
points on the spectral element with acknowledgment that a modeling `crime' has been committed or 
(2) interpolate between the grid (quadrature points) on which the dynamics is calculated to a grid on 
which the physics is calculated (and back), and hence incur an interpolation error.  Since there is a long-standing
history of using the validated physics routines designed for use on uniformly-spaced grids, there is a strong
incentive to apply the second option. However, interpolating density or other key physical quantities without 
accounting for property preservation may lead to negative values that are nonphysical and result in inaccurate representations 
and/or interpretations of the physical data.  For example, Skamrock et al. \cite{skamrock} demonstrated that not preserving positivity may lead to a 
positive bias in a predicted physical quantity of interest (e.g., prediction of moisture).  
The second option mentioned above of moving information from nonuniform to uniform and back via ENO-type interpolation schemes,
explored in \cite{tajo20222PPIsoftware} in the context of high-order methods for numerical
weather prediction, is the main motivation for this work.

Property-preserving interpolation is straightforward when used in the context of low-order numerical simulation methods.
High-order property-preserving interpolation is, however, nontrivial, especially when the interpolation points are not uniformly-spaced. 
In this paper, we demonstrate that it is possible to adaptively construct high-order interpolation methods over unevenly-spaced tensor product grids 
in a way that ensures either data boundedness or positivity preservation (within user-supplied bounds).  
The algorithm we have developed comes with theoretical estimates, presented herein, that provide sufficient conditions for data boundedness and positivity preservation.

\subsection{Previous Work}
\label{sec:prevwork}
In this section, we provide an overview of various numerical approaches to data bounded and positivity preservation.  This overview
is not meant to be exhaustive, but instead to summarize the various ways by which researchers have attempted to tackle this challenging
problem.

Introduced by Harten et al. \cite{HARTEN19973}, Essentially Non-Oscillatory (ENO) schemes were developed to solve problems 
with sharp gradients and discontinuities while achieving high-order accuracy in both smooth and non-smooth regions. As with many 
finite-difference based methods, the backbone of these schemes is interpolation methods.  In the context of this paper, which is to
propose ENO-like interpolation schemes that are property preserving, we briefly review ENO methods.

In the context of finite volume schemes, Fjordholm et al. \cite{fjordholm2013eno} demonstrated that ENO schemes are stable, 
in the sense that the jump of the reconstructed value at each cell interface has the same sign as the jump in the underlying cell average.
Building on the work \cite{tadmor_2003} and \cite{fjordholm2013eno}, Fjordholm et al. \cite{doi:10.1137/110836961} developed a high-order entropy stable ENO scheme for conservation laws.
This approach consists of using entropy conservative flux based on \cite{tadmor_2003}, adding a numerical diffusion to obtain a stable scheme, and obtaining the high-order accuracy via ENO reconstruction.

Harten \cite{harten1989eno,harten1995multiresolution} developed an ENO scheme for subcell resolution in the cases where a discontinuity lies inside a given cell.
Weighted Essentially Non-Oscillatory (WENO) schemes were later proposed by Liu et al. \cite{LIU1994200} to address some of the shortcomings of the ENO schemes.
Shu  \cite{shu_2020} provided a comprehensive overview of different applications and problems in which ENO and WENO schemes are used.  
Shen et al. \cite{SHEN20113780} proposed an adaptive mesh refinement method (AMR) based on WENO schemes for hyperbolic conservation laws. 
In this approach, high-order WENO interpolation is used for the prolongation.
A generalization of the AMR-WENO in \cite{WANG2015161} was used to solve a multi-dimension detonation problem.

Another body of literature sometimes considered around property-preserving methods is computer-aided design and visualization.
Although different from the finite difference (stencil) methods that we seek, we briefly review this literature.
In that literature, ``shape" preservation is often used to describe the preservation of properties like monotonicity and convexity, 
and may include positivity and data boundedness \cite{10.1145/264029.264050} and \cite{Costantini1990}. 
We only briefly review this literature as the additional smoothness constraints at the stencil points enforced by these
methods introduce a level of complexity not needed for our application domain.  Our focus is finite difference ENO-type schemes.
Perhaps the most widely used approach for preserving monotonicity in many applications is PCHIP 
by Fritch and Carlson \cite{fritsch1980monotone}, who derived necessary and sufficient conditions for monotone cubic interpolation, 
and provided an algorithm for building a piecewise cubic approximation from data.
This algorithm calculates the values of the first derivatives at the nodes based on the necessary and sufficient conditions.
Lux et al. \cite{Lux2019ANAF} proposed a monotone quintic spline (MQS) algorithm that relies on the results of He{\ss} and Schmidt \cite{10.1007/BF01934097} and Ulrich and Watson \cite{doi:10.1137/0915035}.	
This method is dependent on the value of the first and second derivatives at the node.
The algorithm uses the sufficient conditions from \cite{10.1007/BF01934097} to check for monotonicity. 
When the conditions are not met, the method in \cite{doi:10.1137/0915035} is used to modify the values of the first and second derivatives to ensure monotonicity.
The work of Dougherty et al. \cite{Hyman} extends these ideas to preserving convexity and  concavity and also to quintic splines.

A second area in which one often finds the development of methods for property preservation is numerical methods for partial differential equations (PDEs). 
Various methods have been developed to enable, for example, positivity-preserving approximations.
To preserve positivity in discontinuous Galerkin (dG) schemes, Zhang et al. \cite{ZHANG2017,Zhang2012,Zhang2752} introduced a linear 
rescaling of polynomials that ensures that the evaluation of the polynomial at the quadrature points remains positive.  
In addition, the linear rescaling of the polynomial conserves mass. 
Light et al. \cite{light} developed a similar approach with a more involved linear polynomial rescaling that preserves positivity at the quadrature nodes and conserves mass.
The polynomial rescaling does not address the case of interpolating between different meshes, which is the primary focus of this work.
Harten et al. \cite{HARTEN19973} developed an Essentially Non-Oscillatory (ENO) piecewise polynomial reconstruction that enables interpolation between different meshes.
The ENO method adaptively chooses stencil points for the interpolation and helps remove Gibbs-like effects but does not guarantee positivity.
As previously mentioned, extensions of these ideas to a Weighted ENO (WENO) combination of these schemes have been proposed by Zhang et al. \cite{Zhang2012_2} and others. 
Finally, Zala et al. \cite{Zala2020,zala2021structurepreserving} developed a nonlinear filtering operator for property-preservation by casting it as an 
optimization problem in which the desired ``structures" (properties) are encoded as constraints.

The  data-bounded interpolation (DBI) method of Berzins \cite{Berzins} builds on three ideas from these ENO and WENO algorithms in the area of the numerical solution of advection equations:
adaptively selecting stencils as in the ENO methods to reduce oscillations \cite{HARTEN19973};
altering the polynomial approximation so that any discontinuities in higher derivatives are removed \cite{harten1995multiresolution}; and
altering the polynomial degree and/or terms so that the ratio of successive divided differences in the series is strictly limited to 
enforce the boundedness of the interpolation \cite{Berzins}.
The work in \cite{berzins2010nonlinear} extends the earlier proof to 1D unevenly-spaced points where, in addition to the interval points, all the remaining points used to build the interpolant are to the right or left of the interval of interest. 
In addition, the work in \cite{berzins2010nonlinear} recognizes that switching off data boundedness when extrema are present is important for maintaining accuracy.
Positivity is important in interpolation cases in which extrema lie between data points and where
the data-bounded interpolant will ``clip" the function, resulting in a loss of accuracy. 
A novel feature of the approach addresses the fact that
preserving positivity alone may still produce undesirable oscillations that lead to an inaccurate representation and/or interpretation of the underlying data.
These oscillations are removed here by imposing strict user-supplied bounds on the positive interpolants as a way of limiting oscillations and correspondingly improving accuracy.

This work extends the ideas in \cite{Berzins} by addressing data boundedness and positivity (within user-supplied bounds) in the same framework and by
allowing meshes of unevenly-spaced points.
The DBI method presented in this paper introduces more relaxed conditions for data-boundedness which give greater accuracy than the conditions used in ~\cite{Berzins}.
Thus, these new proofs provide the previously missing theoretical underpinning for complex interpolation cases such as those
like the NWP case described above. 
The new approach used here both generalizes the
DBI method to unevenly-spaced structured meshes and extends the approach to preserve positivity (positivity-preserving interpolation (PPI)) 
rather than the more restrictive data-bounded approach in \cite{Berzins} and \cite{berzins2010nonlinear}.

\subsection{Outline}
The paper proceeds as follows. In Section \ref{sec:background}, we provide a review of Newton interpolation with particular emphasis on the
properties on Newton polynomials required in this work.  In Section \ref{sec:DBI}, we present our {\em first major contribution}:
theoretical guarantees for adaptive high-order data bounded polynomial interpolation on nonuniformly-spaced points.
In Section \ref{sec:PPI}, we extend the ideas presented in Section \ref{sec:DBI} to positivity.  
We present our {\em second major contribution}:
theoretical guarantees for adaptive high-order positivity-preserving polynomial interpolation on nonuniformly-spaced points.
In Section \ref{sec:PPI}, we also address the case of hidden extrema with the new limiting approach and provided an algorithm for the DBI and PPI methods.
In Section \ref{sec:results}, we provide 1D and 2D results demonstrating the properties of our proposed algorithms.  We summarize and conclude the paper in Section \ref{sec:conclusions}.

\section{Background}
\label{sec:background}

  The approach introduced in this work relies on the Newton polynomial \cite{doi:10.1137/0912034,10.2307/2004888} representation to build interpolants that are positive or bounded by the data values. 
  The ability to adaptively select the divided differences or the stencil as in ENO methods \cite{HARTEN19973} is central to the data-bounded and positivity-preserving interpolation approaches presented in this paper.

  Consider a 1D mesh defined as follows: 
  \begin{equation}\label{eq:RandomMesh}
     \mathcal{M} = \{ x_{i-J}, \cdots, x_{i}, x_{i+1}, \cdots, x_{i+L}\},
  \end{equation} 
  where $x_{i-J} < \cdots < x_{i} < x_{i+1}< \cdots < x_{i+L}$, and $\{ u_{i-J}, \cdots, u_{i+L} \}$ is the set of data values associated with the mesh points.
  In the definition of the mesh $\mathcal{M}$, the subscripts $J$, $L$, $i, \in \mathbb{N}_{0}= \mathbb{N} \cup \{0 \}$, and $x_{k}$, $u_{k} \in \mathbb{R}$ for $i-J \leq k \leq i+L$. 
  For the given mesh $\mathcal{M}$, the Newton divided differences are recursively defined as follows:
  \begin{equation}
    \begin{cases}
    U[x_{i}] = u_{i}\\
    U[x_{i}, \cdots, x_{i+j}] = \frac{U[x_{i+1}, \cdots,x_{i+j}]-U[x_{i}, \cdots,x_{i+j-1}]}{x_{i+j} -x_{i}}.
    \end{cases}
  \end{equation} 

  The ENO procedure starts by setting the initial stencil $\mathcal{V}_{0}$:
  \begin{equation}
    \mathcal{V}_{0} = \{ x_{i}, x_{i+1}\} =\{ x_{0}^{l}, x_{0}^{r}\}. 
  \end{equation}
  The stencil $\mathcal{V}_{0}$ is expanded by successively appending a point to right or left of $\mathcal{V}_{j}$ to form $\mathcal{V}_{j+1}$. 
  The point appended is selected by picking the smallest divided difference at each step.

  Given $\mathcal{V}_{j}$, let $x_{j}^{l}$ and $x_{j}^{r}$ be the leftmost and rightmost stencil points, respectively. 
  In addition, let $x_{p}$ and $x_{q}$ be the stencil points immediately to the left and right of $\mathcal{V}_{j}$. 
  The stencil is expanded from $\mathcal{V}_{j}$ to $\mathcal{V}_{j+1}$ based on the following rules: 
  \begin{itemize}
     \item if $|U[x_{p}, x_{j}^{l}, \cdots, x_{j}^{r}]|<|U[x_{j}^{l}, \cdots, x_{j}^{r}, x_{q}]|$ then \\ 
           $\mathcal{V}_{j+1} =\{ x_{p}, \mathcal{V}_{j} \}$ with $x_{j+1}^{l}=x_{p}$ and $x_{j+1}^{r}=x_{j}^{r}$.
      \item otherwise \\ 
           $\mathcal{V}_{j+1} =\{ \mathcal{V}_{j}, x_{q} \}$ with $x_{j+1}^{l} = x_{j}^{l}$ and $x_{j+1}^{r}= x_{q}$.
  \end{itemize}
  Let 
  \begin{equation}
    I_{i} = [x_{i}, x_{i+1}], \quad \textrm{for } 0 \leq i \leq n-1.
  \end{equation}
  Once the final stencil $\mathcal{V}_{n-1}$ is obtained, the interpolant of degree $n$ defined on $I_{i}$ can be written as
  \begin{equation}\label{eq:newtonPoly}
    \begin{gathered}
    U_{n}(x)= \quad u_{i} + U[x_{0}^{l}, x_{0}^{r}] \pi_{0,i}(x) + 
          U[x_{1}^{l}, \cdots, x_{1}^{r}] \pi_{1,i}(x) + \cdots \\  
          \cdots +  U[x_{n -1}^{l}, \cdots, x_{n -1}^{r}]\pi_{n-1,i}(x),
    \end{gathered}
  \end{equation}
  where 
  $ \pi_{0,i}(x) = (x-x_{i}), \pi_{1,i}(x)=(x-x_{i})(x-x_{1}^{e}), \cdots$ are the Newton basis functions.
  $x_{j}^{e}$ is the point added to expand the stencil $\mathcal{V}_{j-2}$ to $\mathcal{V}_{j-1}$ and can be explicitly expressed as  
  \begin{equation}\label{eq:xe}
     \begin{cases}
       x_{0}^{e} = x_{i}, \\%\quad
       x_{1}^{e} = x_{i+1}, \\% \quad
       x_{j}^{e} = \mathcal{V}_{j-1} \setminus \mathcal{V}_{j-2}, \quad 2 \leq j \leq n-1.
     \end{cases}
  \end{equation}

  The first step in developing the  DBI and PPI methods consists of reorganizing the terms in the polynomial $U_{n}(x)$ defined in Equation (\ref{eq:newtonPoly}) to expose the features used to enforce data boundedness and positivity. 
  The reorganization begins by defining $\lambda_{j}$ as follows:
  \begin{equation}\label{eq:lambda}
    \lambda_{j} =
    \begin{cases}
     1, \quad j = 0 \\
     \frac{U[x_{j}^{l}, \cdots, x_{j}^{r}]}{U[x_{j-1}^{l}, \cdots, x_{j-1}^{r}]}(x_{j}^{r}-x_{j}^{l}), \quad 1\leq j \leq n-1.
    \end{cases}
  \end{equation}
  Expressing $U_{n}(x)$ in terms of $\lambda_{j}$, for $j> 0$ gives
  \begin{equation}\label{eq:UN2}
    \begin{gathered}
    U_{n}(x) = u_{i} + (u_{i+1}-u_{i})\frac{x-x_{0}^{e}}{x_{0}^{r}-x_{0}^{l}} 
              \bigg(1 + \frac{(x-x_{1}^{e})}{(x_{1}^{r}-x_{1}^{l})}\lambda_{1}  \\
              \bigg(1 + \frac{(x-x_{2}^{e})}{(x_{2}^{r}-x_{2}^{l})}\lambda_{2} \bigg( \cdots \lambda_{n-2} 
              \bigg(1 + \frac{(x-x_{n-1}^{e})}{(x_{n-1}^{r}-x_{n-1}^{l})}\lambda_{n-1} \bigg) \cdots \bigg). 
    \end{gathered}
  \end{equation}
  For $x \in I_{i}$, $s$, $t_{j}$, and $d_{j}$ are defined as follows:
  \begin{equation}\label{eq:s}
    0 \leq s = \frac{x-x_{i}}{x_{i+1}-x_{i}} = \frac{x-x_{0}^{e}}{x_{0}^{r}-x_{0}^{l}} \leq 1,
  \end{equation}
  \begin{equation}\label{eq:tk}
    t_{j} = -\frac{x_{i}-x_{j}^{e}}{x_{0}^{r}-x_{0}^{l}}, \textrm{ and } 
  \end{equation}
  \begin{equation}\label{eq:dk}
    0 \leq d_{j} = \frac{x_{j}^{r}-x_{j}^{l}}{x_{0}^{r}-x_{0}^{l}} .
  \end{equation}
  $s$ and $d_{j}$ are defined such that $s \in [0, 1]$ and $d_{j} \geq 0$.
  Expressing $\frac{x-x_{j}^{e}}{x_{j}^{r}-x_{j}^{l}}$ in terms of $s$, $t_{j}$, and $d_{j}$ gives  
  \begin{equation}\label{eq:stk}
    \frac{x-x_{j}^{e}}{x_{j}^{r}-x_{j}^{l}} = \frac{ \frac{x-x_{i}}{x_{0}^{r}-x_{0}^{l}} +  \frac{x_{i}-x_{j}^{e}}{x_{0}^{r}-x_{0}^{l}}}
                                               { \frac{x_{j}^{r}-x_{j}^{l}}{x_{0}^{r}-x_{0}^{l}} }
                                               = \frac{s-t_{j}}{d_{j}}.
  \end{equation}
  Using the results from Equation (\ref{eq:stk}), the polynomial $U_{n}(x)$ as expressed in Equation (\ref{eq:UN2}) can be written as
  \begin{equation}\label{eq:UN3}
    \begin{gathered}
    U_{n}(x) = u_{i} + (u_{i+1}-u_{i})S_{n}(x)
    \end{gathered}
  \end{equation} 
  with $S_{n}(x)$ defined as
  \begin{equation}\label{eq:SN1}
    \begin{gathered}
    S_{n}(x) = s
              \bigg(1 + \frac{(s-1)}{d_{1}}\lambda_{1} 
              \bigg(1 + \frac{(s-t_{2})}{d_{2}}\lambda_{2} \bigg( \cdots %\lambda_{n-2} 
              \bigg(1 + \frac{(s-t_{n-1})}{d_{n-1}}\lambda_{n-1} \bigg) \cdots \bigg). 
    \end{gathered}
  \end{equation} 
  For future use below, $S_{n}(x)$ can be compactly represented by introducing $\delta_{j}$ defined as
  \begin{equation} \label{eq:deltaj}
    \begin{cases}
      \delta_{n} = 1\\
      \delta_{j} = 1+ \frac{s-t_{j}}{d_{j}}\lambda_{j}\delta_{j+1} \quad 2\leq j \leq n-1 \\
      \delta_{1} = s+ \frac{s(s-1)}{d_{1}}\delta_{2} = S_{n}(x).
    \end{cases}
  \end{equation}
  Together, $U_{n}(x)$ and $S_{n}(x)$ in Equations (\ref{eq:UN3}) and (\ref{eq:SN1}) are reorganizations needed to construct the DBI and PPI algorithm.
  The general approach is to first bound the quadratic term in $S_{n}(x)$ and then to increase the order to cubic, quartic, and higher order polynomials.
  This iterative procedure is used to define computational bounds on the values of $\bar{\lambda}_{j} = \prod_{k=0}^{j}\lambda_{k}$.
  $\bar{\lambda}_{j}$ can be explicitly written  as
  \begin{equation} \label{eq:bar_lambda}
    \bar{\lambda}_{j} = \lambda_j \bar{\lambda}_{j-1} = \prod_{k=1}^{j}\lambda_{k} =
    \begin{cases}
      1 \quad j = 0,\\
    \frac{U[x_{j}^{l}, \cdots, x_{j}^{r}]}{U[x_{0}^{l}, x_{0}^{r}]}\prod_{k=1}^{j}(x_{k}^{r}-x_{k}^{l}), \quad 1\leq j \leq n-1.
    \end{cases}
  \end{equation}

 \section{Data-Bounded Interpolation}
  \label{sec:DBI}
  The DBI method builds on three ideas from algorithms in the area of the numerical solution of advection equations: 
  adaptively selecting stencils as in the ENO methods to reduce oscillations \cite{HARTEN19973}; 
  altering the polynomial approximation so that any discontinuities in higher derivatives are removed \cite{harten1995multiresolution}; and
  altering the polynomial degree and/or terms so that the ratio of successive divided differences in the series is strictly limited to enforce the boundedness of the interpolation \cite{Berzins}.  
  In the DBI method introduced here, more relaxed bounds on $\bar{\lambda_{j}}$ defined in Equation (\ref{eq:bar_lambda}) are derived which gives greater accuracy than those in ~\cite{Berzins}. 
  The work in ~\cite{Berzins} requires that the absolute values of $\bar{\lambda}_{j}$ decrease as more terms are added ($|\bar{\lambda}_{j}|>|\bar{\lambda}_{j+1}|$) and $|\bar{\lambda}_{j}|<1$ which are more restrictive than the bounds in Equation (\ref{eq:bound}).
  For a given set of mesh points and the data values associated with those mesh points, we approximate the data with a $\mathbf{C}^{0}$ continuous function that is built by fitting a polynomial in each subinterval $I_{i}$.
  The fitted polynomial is constructed in such a way that it is bounded by $u_{i}$ and $u_{i+1}$. 
  Given that this work concerns itself with locally fitting a polynomial in the interval $I_{i}$, let us assume, for the remaining parts of this paper, that $x\in I_{i}$ and that building the interpolant always starts with the stencil $\mathcal{V}_{0} = \{ x_{i}, x_{i+1}\}$.
 
  Let $U^{l}(x)$ be the limited polynomial defined as in Equation (\ref{eq:UN3}) and bounded by $u_{i}$ and $u_{i+1}$.
  For the polynomial $U^{l}(x)$ to be bounded by $u_{i}$ and $u_{i+1}$, it follows that  for $x \in I_{i}$
  \begin{equation}\label{eq:SNBounds}
     0 \leq S_{n}(x) \leq 1,
  \end{equation}
  with $S_{n}(x)$ defined in Equation (\ref{eq:SN1}).
  The reconstruction procedure begins by considering the linear and quadratic terms from $S_{n}(x)$ in Equation (\ref{eq:SN1}), and imposing the following bounds:
  \begin{equation}\label{eq:quad1}
    0 \leq s\big(1 + \frac{s-1}{d_{1}} \bar{\lambda}_{1}\big) \leq 1. 
  \end{equation}
  As $s \in [0,1]$ and isolating $\bar{\lambda}_{1}$ in Equation (\ref{eq:quad1}) gives 
  \begin{equation}
      -\frac{d_{1}}{s} \leq \bar{\lambda}_{1} \leq \frac{d_{1}}{1-s}, \textrm{ and } 
  \end{equation}
  \begin{equation} \label{eq:lambda_bar_1}
    -d_{1} \leq \bar{\lambda}_{1} \leq d_{1}.
  \end{equation}
  The bounds from Equation (\ref{eq:lambda_bar_1}) are extended to bound the cubic form by requiring 
that what multiplies $\bar{\lambda}_{1}$ must fit into the inequality in Equation (\ref{eq:lambda_bar_1}).  
  Thus, for the cubic case Equation (\ref{eq:lambda_bar_1}) becomes 
  \begin{equation} \label{eq:lambda_bar_2_2}
    -d_{1} \leq \bar{\lambda}_{1}\big(1 + \frac{(s-t_{2})}{d_{2}}\lambda_{2} \big)\leq d_{1}.
  \end{equation}
  Subtracting $\bar{\lambda}_{1}$ from this inequality gives 
  \begin{equation}
    -d_{1} -\bar{\lambda}_{1}\leq \frac{(s-t_{2})}{d_{2}}\bar{\lambda}_{2} \leq d_{1}-\bar{\lambda}_{1}.
  \end{equation}
  In the case when $t_{2}$ is negative, $s-t_{2}$ has a maximum value at $s=1$ and a minimum value at $s=0$. 
  $\bar{\lambda}_{2}$ is then bounded by 
  \begin{equation}\label{eq:lambda_bar_2}
    \frac{d_{2}}{(1-t_{2})}(-d_{1} -\bar{\lambda}_{1}\big)\leq\bar{\lambda}_{2}\leq (d_{1}-\bar{\lambda}_{1}) \frac{d_{2}}{(1-t_{2})}.
  \end{equation}
  When $t_{2}$ positive, $\frac{1}{1-t_{2}}$ is substituted by $\frac{1}{-t_{2}}$ and the inequalities 
  $\leq$ with $\geq$ and vice versa are swapped. 
  In the quartic case, we require that 
  \begin{equation}
     \frac{d_{2}}{1-t_{2}}(-d_{1} -\bar{\lambda}_{1})
     \leq \bar{\lambda}_{2}\bigg(1 +\frac{(s-t_{3})}{d_{3}}\lambda_{3} \bigg)
     \leq \frac{d_{2}}{1-t_{2}} (d_{1}-\bar{\lambda}_{1}). 
  \end{equation}
  If we assume that $t_{3}$ is negative
  \begin{equation}\label{eq:lambda_bar_3}
    \frac{d_{3}}{1-t_{3}} \bigg(\frac{d_{2}}{1-t_{2}}(-d_{1} -\bar{\lambda}_{1}\big) - \bar{\lambda}_{2}\bigg)
    \leq \bar{\lambda}_{3} \leq 
    \frac{d_{3}}{1-t_{3}}\bigg( \frac{d_{2}}{1-t_{2}} (d_{1}-\bar{\lambda}_{1}) -\bar{\lambda}_{2} \bigg). 
  \end{equation}
  This reconstruction procedure can be continued to higher orders provided that care is taken to correctly manage the impact of the signs of $t_{j}$. 
  For the boundary and nearby boundary intervals, fewer choices are available, and the final stencil is biased towards the interior of the domain because there are no points to choose from beyond the boundaries. 
  In the process of constructing $\mathcal{V}_{n-1}$, when the left or right boundary are reached, the remaining mesh points are obtained from the side that is towards the interior of the domain.

  For a more formal  and complete expression of this recursive procedure, the bounds on $\bar{\lambda}_{j}$ can be defined as follows:
  \begin{subequations}
    \begin{equation}\label{eq:Bminus}
        B_{j}^{-} = 
      \begin{cases}
        -d_{1} \quad  & j =0 \\
        (B_{j-1}^{-} -\bar{\lambda}_{j-1}) \frac{d_{j}}{1-t_{j}}, & t_{j} \in (-\infty, 0] \quad j>1 \\ 
        (B_{j-1}^{+} -\bar{\lambda}_{j-1}) \frac{d_{j}}{-t_{j}},  & t_{j} \in (0, +\infty) \quad j>1, 
      \end{cases}
    \end{equation}
    \textrm{and}
    \begin{equation}\label{eq:Bplus}
      B_{j}^{+} =
      \begin{cases}
         d_{1} , &  j = 1\\
         (B_{j-1}^{+} -\bar{\lambda}_{j-1}) \frac{d_{j}}{1-t_{j}}, & t_{j} \in (-\infty, 0] \quad j >1\\ 
         (B_{j-1}^{-} -\bar{\lambda}_{j-1}) \frac{d_{j}}{-t_{j}}, &  t_{j} \in (0, +\infty) \quad j >1.
      \end{cases}
    \end{equation}
 \end{subequations}
  The sign of $t_{j}$ is incorporated into the definitions of $B_{j}^{-}$ and $B_{j}^{+}$ in Equations (\ref{eq:Bminus}) and (\ref{eq:Bplus}), respectively.
  The sufficient conditions for data boundedness such as 
  Equations (\ref{eq:lambda_bar_1}), (\ref{eq:lambda_bar_2}) and  (\ref{eq:lambda_bar_3})  can now be written as
  \begin{equation}\label{eq:bound}
    B_{j}^{-} \leq \bar{\lambda}_{j} \leq B_{j}^{+}, \textrm{ for } j \geq 0.
  \end{equation}

  %%%%%%%%%%%%%%%%%%%%%%%%%%%%%%%%%%%%%%%%%%%%%%%%%%%%%%%%%%%%%%%%%%%%%%%%%%%%%%%%%%%%%%%%%%%%%%%%%%%%%%%%%%%%%%%%%%%%%%%%%%%%%%%%%%%%%%%%%%%%%%%%%%
  \begin{lemma}\label{lem:DBI}
    Let us assume that for $x \in I_{i}$, $B_{j}^{-}$ and $B_{j}^{+}$ are defined as in Equations (\ref{eq:Bplus}) and (\ref{eq:Bminus}), respectively.
    In addition, let $\delta_{j}$ be defined as in Equation (\ref{eq:deltaj}).
    If for $x \in I_{i}$, $B_{j}^{-}$ is negative, $B_{j}^{+}$ is positive, and $B_{j}^{-} \leq \bar{\lambda}_{j}\delta_{j+1} \leq B_{j}^{+}$, then 
    $$B_{j-1}^{-} \leq \bar{\lambda}_{j-1}\delta_{j} \leq B_{j-1}^{+}.$$
  \end{lemma}
  \begin{proof}
    The proof is split into two cases that take into consideration the different possible values of $t_{j}$, and in each case we consider the left and right side of the inequality separately.
    \begin{enumerate}[label=\textbf{(\Roman*)}]
      \item $t_{j} \in (-\infty , 0]$ \\
        Let us start with the left side of the inequality (i.e., $B_{j-1}^{-} \leq \bar{\lambda}_{j-1}\delta_{j}$). 
        Noting that $\frac{1-t_{j}}{s-t_{j}} \geq 1$ for $s \in [0,1]$, and using $B_{j}^{-} \leq 0 $ and $B_{j}^{-} \leq \bar{\lambda}_{j}\delta_{j+1}$, we have
        \begin{equation}\label{eq:Bj_left1}
          \begin{aligned}
            (B_{j-1}^{-} -\bar{\lambda}_{j-1})\frac{d_{j}}{s-t_{j}}  =& \frac{1-t_{j}}{s-t_{j}} B_{j}^{-}\\
                                                                     \leq    & B_{j}^{-} \\ 
                                                                     \leq & \bar{\lambda}_{j}\delta_{j+1}. 
          \end{aligned}
        \end{equation}
        Isolating $B_{j-1}^{-}$ in Equation (\ref{eq:Bj_left1}) and using Equations (\ref{eq:deltaj}) and (\ref{eq:bar_lambda}) leads to 
        \begin{equation}\label{eq:Bj-1_left1}
          \begin{aligned}
            B_{j-1}^{-} \leq & \bar{\lambda}_{j-1}+\frac{s-t_{j}}{d_{j}}\bar{\lambda}_{j}\delta_{j+1} \\
                        \leq & \bar{\lambda}_{j-1}\bigg(1+\frac{s-t_{j}}{d_{j}}\lambda_{j}\delta_{j+1} \bigg) \\
                        =    & \bar{\lambda}_{j-1}\delta_{j}.
          \end{aligned}
        \end{equation}
        Now, let us focus on the right side of the inequality (i.e., $B_{j-1}^{+} \geq \bar{\lambda}_{j-1}\delta_{j}$)
        Again, observing that $\frac{1-t_{j}}{s-t_{j}} \geq 1$ for $s \in [0,1 ]$ and using $B_{j}^{+} \geq 0$ and $B_{j}^{+} \geq \bar{\lambda}_{j}\delta_{j+1}$ yields
        %%%
        \begin{equation}\label{eq:Bj_right1}
          \begin{aligned}
            (B_{j-1}^{+}-\bar{\lambda}_{j-1})\frac{d_{j}}{s-t_{j}}      = & \frac{1-t_{j}}{s-t_{j}}B_{j}^{+} \\
                                              	                     \geq & B_{j}^{+} \\
                                                                     \geq & \bar{\lambda}_{j}\delta_{j+1}.
          \end{aligned}
        \end{equation}
        Isolating $B_{j-1}^{+}$ in Equation (\ref{eq:Bj_right1}) yields 
        \begin{equation}\label{eq:Bj-1_right1}
          \begin{aligned}
            B_{j-1}^{+} \geq & \bar{\lambda}_{j-1}+\frac{s-t_{j}}{d_{j}}\bar{\lambda}_{j}\delta_{j+1} \\
                        \geq & \bar{\lambda}_{j-1}\bigg(1+\frac{s-t_{j}}{d_{j}}\lambda_{j}\delta_{j+1} \bigg) \\
                          =  & \bar{\lambda}_{j-1}\delta_{j}.
          \end{aligned}
        \end{equation}
      \item $t_{j} \in (0, + \infty)$

        Let us consider the left side of the inequality (i.e., $B_{j-1}^{-} \leq \bar{\lambda}_{j-1}\delta_{j}$). 
        Multiplying $B_{j}^{-}$ by $\frac{-t_{j}}{s-t_{j}}$ yields
        \begin{equation}\label{eq:Bj_left}
          \begin{aligned}
            (B_{j-1}^{+} -\bar{\lambda}_{j-1})\frac{d_{j}}{s-t_{j}} = & \frac{-t_j}{s-t_j} B_{j}^{-}. 
          \end{aligned}
        \end{equation}
        Given that $B_{j}^{-} \leq 0 $ and $B_{j}^{-} \leq \bar{\lambda}_{j}\delta_{j+1}$, and noting that $\frac{-t_{j}}{s-t_{j}} \geq 1$ for  $s \in [0, 1]$, the right side of Equation (\ref{eq:Bj_left}) can be bounded by $B_{j}^{-}$ to give 
        \begin{equation}\label{eq:Bj_left2}
          \begin{aligned}
            (B_{j-1}^{+} -\bar{\lambda}_{j-1})\frac{d_{j}}{s-t_{j}} \leq & B_{j}^{-} \\
                                                                    \leq & \bar{\lambda}_{j}\delta_{j+1}.
          \end{aligned}
        \end{equation}
        Isolating $B_{j-1}^{+}$ in Equation (\ref{eq:Bj_left2}) leads to 
        \begin{equation}\label{eq:Bj-1_left2}
          \begin{aligned}
          B_{j-1}^{+} \geq & \bar{\lambda}_{j-1}+\frac{s-t_{j}}{d_{j}}\bar{\lambda}_{j}\delta_{j+1} \\
                      \geq & \bar{\lambda}_{j-1}\bigg(1+\frac{s-t_{j}}{d_{j}}\lambda_{j}\delta_{j+1} \bigg)\\
                         = & \bar{\lambda}_{j-1}\delta_{j}.
          \end{aligned}
        \end{equation}
        For the right side of the inequality (i.e. $B_{j-1}^{-} \leq \bar{\lambda}_{j-1}\delta_{j}$), 
        $\frac{-t_{j}}{s-t_{j}} \geq 1$ for $s \in [0,1]$, and using $B_{j}^{-} \leq 0 $ and $B_{j}^{-} \leq \bar{\lambda}_{j}\delta_{j+1}$ yields 
        \begin{equation}\label{eq:Bj_right2}
          \begin{aligned}
            (B_{j-1}^{-} -\bar{\lambda}_{j-1})\frac{d_{j}}{s-t_{j}} = & \frac{-t_j}{s-t_j} B_{j}^{+} \\
                                                                 \geq & B_{j}^{-} \\
                                                                 \geq & \bar{\lambda}_{j}\delta_{j+1}.
          \end{aligned}
        \end{equation}
        Isolating $B_{j-1}^{-}$ in Equation (\ref{eq:Bj_right2}) yields
        \begin{equation}\label{eq:Bj-1_right2}
          \begin{aligned}
          B_{j-1}^{-}  \leq & \bar{\lambda}_{j-1}+\frac{s-t_{j}}{d_{j}}\bar{\lambda}_{j}\delta_{j+1} \\
                       \leq & \bar{\lambda}_{j-1}\bigg(1+\frac{s-t_{j}}{d_{j}}\lambda_{j}\delta_{j+1} \bigg) \\
                          = & \bar{\lambda}_{j-1}\delta_{j}.
          \end{aligned}
        \end{equation}
      \end{enumerate}
      The results from Equations (\ref{eq:Bj-1_left1}), (\ref{eq:Bj-1_right1}), (\ref{eq:Bj-1_left1}), and (\ref{eq:Bj-1_right1}) can be summarized as 
      \begin{equation*}
        B_{j-1}^{-} \leq \bar{\lambda}_{j-1}\delta_{j} \leq B_{j-1}^{+}.
      \end{equation*}
  \end{proof}
  %%%%%%%%%%%%%%%%%%%%%%%%%%%%%%%%%%%%%%%%%%%%%%%%%%%%%%%%%%%%%%%%%%%%%%%%%%%%%%%%%%%%%%%%%%%%%%%%%%%%%%%%%%%%%%%%%%%%%%%%%%%%%%%%%%%%%%%%%%%%%%%%%%
  \begin{theorem}\label{theo:DBI}
    Assuming that for $x \in I_{i}$, the polynomial $S_{n}(x)$ of degree $n$ is built starting from the stencil $\mathcal{V}_{0} = \{x_{i}, x_{i+1}\}$,
and then by successively appending mesh points from the left and/or right of the interval $I_{i}$ to obtain the final stencil $\mathcal{V}_{n-1}$. 
The construction of $\mathcal{V}_{n-1}$ does not require the points to be added in a symmetric fashion alternating from left to right.
    If for $x \in I_{i}$, $B_{j}^{-}$ defined in Equation (\ref{eq:Bminus}) is negative, $B_{j}^{+}$ defined in Equation (\ref{eq:Bplus}) is positive, and $B_{j}^{-} \leq \bar{\lambda}_{j} \leq B_{j}^{+}$ then for $x \in I_{i}$ 
    $$0 \leq S_{N}(x) \leq 1.$$
  \end{theorem}
  \begin{proof}
    This proof builds on the results from Lemma \ref{lem:DBI} and  
    starts by using $B_{j}^{-} \leq \bar{\lambda}_{j} \leq B_{j}^{+}$ to bound $\bar{\lambda}_{n-1}$ as follows:
    \begin{equation}\label{eq:theo_lambda_n-1}
      B_{n-1}^{-} \leq \bar{\lambda}_{n-1} \leq B_{n-1}^{+}.
    \end{equation}
    By Lemma \ref{lem:DBI}, Equation (\ref{eq:theo_lambda_n-1}) then leads to
    \begin{equation} 
      B_{n-2}^{-} \leq \bar{\lambda}_{n-2}\delta_{n-1} \leq B_{n-2}^{+}. 
    \end{equation}
    Successively, using the results from Lemma \ref{lem:DBI} to bound $\bar{\lambda}_{n-2}\delta_{n-1}$, $\bar{\lambda}_{n-3}\delta_{n-2}$, $\cdots$, $\bar{\lambda}_{1}\delta_{2}$, yields  
    \begin{equation} \label{eq:theo_lambda_1}
      B_{1}^{-} \leq \bar{\lambda}_{1}\delta_{2} \leq B_{1}^{+}, 
    \end{equation} 
    where $\delta_{j}$ is defined in Equation (\ref{eq:deltaj}).
    The results from Equation (\ref{eq:theo_lambda_1}) may now be used to derive the target bounds (i.e., $0 \leq S_{N}(x) \leq 1$).
    Considering the left side of Equation (\ref{eq:theo_lambda_1}) (i.e., $B_{1}^{-} \leq \bar{\lambda}_{1}\delta_{2}$), and noting that $\frac{(s-1)}{s(s-1)} \geq 1$, gives 
    \begin{equation}\label{eq:delta2_left}
      \begin{aligned}
        -\frac{(s-1)}{s(s-1)} d_{1}=& B_{1}^{-}\frac{(s-1)}{s(s-1)} \\
                               \leq & B_{1}^{-}\\
                               \leq & \bar{\lambda}_{1}\delta_{2}. 
      \end{aligned}
    \end{equation}
    Isolating $\delta_{1}$ from Equation (\ref{eq:delta2_left}) gives
    \begin{equation}\label{eq:delta1_left}
      1 \geq s+\frac{s(1-s)}{d_{1}}\bar{\lambda}_{1}\delta_{2} = \delta_{1} = S_{n}(x). 
    \end{equation}
    Considering the right side of Equation (\ref{eq:theo_lambda_1}) (i.e. $B_{1}^{+} \geq \bar{\lambda}_{1}\delta_{2}$), and noting that $\frac{(-s)}{s(s-1)} \geq 1$, gives 
    \begin{equation}\label{eq:delat2_right}
      \begin{aligned}
        \frac{(-s)}{s(s-1)} d_{1}=& B_{1}^{+}\frac{(-s)}{s(s-1)} \\
                               \geq & B_{1}^{+}\\
                               \geq & \bar{\lambda}_{1}\delta_{2}. 
      \end{aligned}
    \end{equation}
    Isolating $\delta_{1}$ from Equation (\ref{eq:delat2_right}) gives
    \begin{equation}\label{eq:delta1_right}
      0 \leq s+\frac{s(1-s)}{d_{1}}\bar{\lambda}_{1}\delta_{2} = \delta_{1} = S_{n}(x). 
    \end{equation}
    The proof concludes by combining the results from Equations (\ref{eq:delta1_left}) and (\ref{eq:delta1_right}) to obtain 
    \begin{equation}
      0 \leq s+\frac{s(1-s)}{d_{1}}\bar{\lambda}_{1}\delta_{2} = \delta_{1} = S_{n}(x) \leq 1.
    \end{equation}
  \end{proof}

  \section{Constrained Positivity-Preserving Interpolation}
  \label{sec:PPI}

  In many cases, it is sufficient to preserve positivity through interpolation and not to enforce the stricter requirement of data boundedness.
  As mentioned in the introduction, the case of unknown extrema between data points is an important example. 
  Let $U^{p}(x)$ be a positive polynomial of degree $n$ defined over the interval $I_{i}$ as in Equation (\ref{eq:UN3}).
  For $x \in I_{i}$, the polynomial $U^{p}(x)$ is allowed to grow beyond $u_{i}$ and $u_{i+1}$ but must remain positive. 
  For the polynomial to be positive, one requires that 
  \begin{equation}\label{eq:pos}
    U^{p}(x) \geq 0.
  \end{equation}
  However, in practice, enforcing positivity alone may still result in large oscillations and in extrema that degrade the approximation.
  We observe this behavior because enforcing positivity alone does  not restrict how much the polynomial is allowed to grow beyond the data values. 
  In addition to enforcing positivity, it is important to remove the undesirable oscillations and extrema as much as possible.
  Let us define $u_{min}$ and $u_{max}$ as
  \begin{equation}\label{eq:umin}
     u_{min} =  \mathbf{min} ( u_{i}, u_{i+1} ) -\Delta_{min}, 
  \end{equation}
  and 
  \begin{equation}\label{eq:umax}
     u_{max} = \mathbf{max} (u_{i}, u_{i+1} ) + \Delta_{max},
  \end{equation}
  where $\Delta_{min}$ and $\Delta_{max}$ are user-defined parameters used to bound the positive polynomial $U^{p}(x)$.
  To allow the polynomial to grow beyond the data values but not produce extrema that are too large, we bound $U^{p}(x)$ as follows:
  \begin{equation}\label{eq:boundUP}
    u_{min} \leq U^{p}(x)= u_{i} + (u_{i+1}-u_{i})S_{n}(x) \leq u_{max}.
  \end{equation} 
   The interpolant $U^{p}(x)$ is now positive and bounded by $u_{min}$ and $u_{max}$.
  Equation (\ref{eq:boundUP}) is equivalent to bounding $S_{n}(x)$ as follows:
  \begin{equation}
    m_{\ell} \leq S_{n}(x) \leq m_r, 
  \end{equation} 
  where the factors $m_{\ell}$ and $m_{r}$ are expressed as 
  \begin{enumerate}[label=\textbf{(\Roman*)}]
    \item : $u_{i+1} > u_{i}$
      \begin{equation}\label{eq:CaseI}
        m_{\ell} = \mathbf{min}\bigg( 0, \frac{u_{min}-u_{i}}{u_{i+1}-u_{i}}\bigg), \textrm{ and }
        m_{r} = \mathbf{max}\bigg( 1, \frac{u_{max}-u_{i}}{u_{i+1}-u_{i}}\bigg)
      \end{equation}
    \item : $u_{i+1} < u_{i}$
      \begin{equation}\label{eq:CaseII}
        m_{\ell} = \mathbf{min}\bigg( 0, \frac{u_{max}-u_{i}}{u_{i+1}-u_{i}}\bigg), \textrm{ and }
        m_{r} = \mathbf{max}\bigg( 1, \frac{u_{min}-u_{i}}{u_{i+1}-u_{i}}\bigg).
      \end{equation}
  \end{enumerate}
  We note that if we set $\Delta_{min} = 0$ and $\Delta_{max}=0$, we recover Equation (\ref{eq:pos}).

  The PPI method is constructed by relaxing the bounds imposed on $\bar{\lambda}_{1}$ as follows:
  \begin{equation}\label{eq:lambda1_bounds_ppi}
    \bigg(-4(m_{r}-1)-1\bigg)d_{1} \leq \bar{\lambda}_{1} \leq \bigg(-4m_{\ell} + 1 \bigg)d_{1}.
  \end{equation}
  Let us demonstrate how the PPI method is constructed in the case of a quadratic interpolant.
  Starting from the DBI results in the Theorem \ref{theo:DBI}, it follows that 
  \begin{equation}\label{eq:ppi_quad}
    0 \leq s + \frac{s(s-1)}{d_{1}}\bar{\lambda}_{1} \leq 1.
  \end{equation}
  Relaxing the left and right bounds in Equation (\ref{eq:ppi_quad}) by $m_{\ell}$ and $m_{r}$, respectively leads to
  \begin{equation}\label{eq:ppi_quad2}
    m_{\ell} \leq s + \frac{s(s-1)}{d_{1}}\bar{\lambda}_{1} \leq m_{r}.
  \end{equation}
  Isolating $\bar{\lambda}_{1}$ from Equation (\ref{eq:ppi_quad2}) leads to
  \begin{equation}\label{eq:ppi_quad3}
    \frac{m_{r}-s}{s(s-1)}d_{1} \leq \bar{\lambda}_{1} \leq \frac{m_{\ell}-s}{s(s-1)}d_{1}.
  \end{equation}
  Equation (\ref{eq:ppi_quad3}) can be reorganized to obtain
  \begin{equation}
    \bigg(\frac{m_{r}-1}{s(s-1)}+\frac{1-s}{s(s-1)}\bigg)d_{1} 
    \leq \bar{\lambda}_{1} \leq 
    \bigg(\frac{m_{\ell}}{s(s-1)} - \frac{s}{s(s-1)} \bigg)d_{1} 
  \end{equation}
  and then %%
  \begin{equation}
    \bigg(\frac{m_{r}-1}{s(s-1)}-\frac{1}{s}\bigg)d_{1} 
    \leq \bar{\lambda}_{1} \leq 
    \bigg(\frac{m_{\ell}}{s(s-1)} - \frac{1}{(s-1)} \bigg)d_{1}.
  \end{equation}
  Noting that $\frac{1}{s(s-1)} \leq -4$, $\frac{1}{s} \geq 1$, and $\frac{1}{s-1} \leq -1$, we obtain
  \begin{equation}\label{eq:lambda_bar_ppi_1}
    \bigg(-4(m_{r}-1)-1\bigg)d_{1} 
    \leq \bar{\lambda}_{1} \leq 
    \bigg(-4m_{\ell} + 1 \bigg)d_{1}.
  \end{equation}
  Once the bounds on $\bar{\lambda}_{1}$ and  the quadratic interpolant are determined, the extension to cubic, quartic, and higher order interpolants follows the same reconstruction procedure used in the DBI method and outlined from Equation (\ref{eq:lambda_bar_2_2}) to (\ref{eq:lambda_bar_3}).
  As in the case of the DBI method, fewer choices are available for $\mathcal{V}_{n-1}$ at the boundary and nearby boundary intervals because there are no points to choose from beyond the boundaries. 
  When a boundary is reached during the process of constructing the stencil $\mathcal{V}_{n-1}$, the remaining mesh points are picked from the side that is towards the interior of the domain.
  The final stencil at the boundary and nearby the boundary intervals are biased towards the interior of the domain. 
  The recursive expression for the bounds on $\bar{\lambda}_{j}$ for the PPI method becomes
   \begin{subequations}
   \begin{equation}\label{eq:B_PPIminus}
        B_{j}^{-} = 
      \begin{cases}
        (-4(m_{r}-1) -1)d_{1} \quad j =1 \\
        (B_{j-1}^{-} -\bar{\lambda}_{j-1}) \frac{d_{j}}{1-t_{j}}, \textrm{ if } t_{j} \in (-\infty, 0] \quad j>1 \\ 
        (B_{j-1}^{+} -\bar{\lambda}_{j-1}) \frac{d_{j}}{-t_{j}}, \textrm{ if } t_{j} \in (0, 1) \cup (1, +\infty) \quad j>1, 
      \end{cases}
    \end{equation}
    \textrm{and}
    \begin{equation}\label{eq:B_PPIplus}
      B_{j}^{+} =
      \begin{cases}
         (-4m_{\ell}+1)d_{1} , \quad j = 1\\
         (B_{j-1}^{+} -\bar{\lambda}_{j-1}) \frac{d_{j}}{1-t_{j}}, \textrm{ if } t_{j} \in (-\infty, 0] \quad j >1\\ 
         (B_{j-1}^{-} -\bar{\lambda}_{j-1}) \frac{d_{j}}{-t_{j}}, \textrm{ if } t_{j} \in (0, +\infty) \quad j >1.
      \end{cases}
    \end{equation}
 \end{subequations}
  The difference between the DBI and PPI methods is highlighted in how the bounds $B_{1}^{-}$ and $B_{1}^{+}$ are calculated.
  More precisely, $B_{1}^{-}$  and $B_{1}^{+}$ are defined as $-d_{1}$ and $d_{1}$ for the DBI method, whereas for the PPI method, they are defined as $(-4(m_{r}-1) -1)d_{1}$ and $(-4m_{\ell}+1)d_{1}$, respectively.
  In addition, the DBI method can be recovered from the PPI methods by setting $m_{\ell} = 0$ and $m_{r} = 1$.
  For example, in the case of the right boundary Equations (\ref{eq:lambda_bar_1}) and (\ref{eq:lambda_bar_ppi_1}) can be written as
  \begin{equation} \label{eq:lambda_bar_dbi_2}
    -d_{1} \leq \bar{\lambda}_{1} = \frac{U[x_{N-2}, x_{N-1}, x_{N}]}{U[x_{N-1}, x_{N}]}(x_{N}-x_{N-1})\leq d_{1}, \textrm{ and }
  \end{equation}
  \begin{equation}\label{eq:lambda_bar_ppi_2}
    \bigg(-4(m_{r}-1)-1\bigg)d_{1} 
    \leq \bar{\lambda}_{1} =\frac{U[x_{N-2}, x_{N-1}, x_{N}]}{U[x_{N-1}, x_{N}]}(x_{N}-x_{N-1})\leq 
    \bigg(-4m_{\ell} + 1 \bigg)d_{1},
  \end{equation}
  where $x_{N}$ is the mesh point at the right boundary, $m_{\ell} \leq 0$, $m_{r} \geq 1$, and  
  \begin{equation}
    d_{1} = \frac{x_{N}-x_{N-2}}{x_{N}-x_{N-1}}.
  \end{equation}
   From Equations (\ref{eq:CaseI}) and (\ref{eq:CaseII}), $m_{r}=18.94$ and $m_{\ell}=-18.94$ for the right boundary of the Runge example in Figure \ref{fig:runge_oscillations} below.
  Equations (\ref{eq:lambda_bar_dbi_2}) and (\ref{eq:lambda_bar_ppi_2}) show the bounds on $\bar{\lambda}_{1}$ for data-boundedness and positivity, respectively.
  Given that $(-4(m_{r}-1)-1) \leq 0$ and $(-4m_{\ell} + 1) \geq 1$, the bounds for positivity are more relaxed than data-boundedness.
  Thus, enabling the use of higher degree polynomials for the PPI method than for the DBI method.
  %%%%%%%%%%%%%%%%%%%%%%%%%%%%%%%%%%%%%%%%%%%%%%%%%%%%%%%%%%%%%%%%%%%%%%%%%%%%%%%%%%%%%%%%%%%%%%%%%%%%%%%%%%%%%%%%%%%%%%%%%%%%%%%
  \begin{theorem}\label{theo:SN}
    Let us assume that for $x \in I_{i}$, the polynomials $U_{n}(x)$ and $S_{n}(x)$ of degree $n$ are defined as in Equations (\ref{eq:UN3}) and (\ref{eq:SN1}), respectively. Both polynomials are built starting from the stencil $\mathcal{V}_{0} = \{x_{i}, x_{i+1}\}$,
and then by successively appending mesh points from the left and/or right of the interval $I_{i}$ to obtain the final stencil $\mathcal{V}_{n-1}$. 
The construction of $\mathcal{V}_{n-1}$ does not require the points to be added in a symmetric fashion alternating from left to right.
    If for $x \in I_{i}$, $B_{j}^{-}$ defined in Equation (\ref{eq:Bminus}) is negative, $B_{j}^{+}$ defined in Equation (\ref{eq:Bplus}) is positive, and $B_{j}^{-} \leq \bar{\lambda}_{j} \leq B_{j}^{+}$ then for $x \in I_{i}$ 
      $$m_{\ell} \leq S_{n}(x) \leq m_r,$$
    where $m_{\ell}$ and $m_r$ are provided in Equations (\ref{eq:CaseI}) and (\ref{eq:CaseII}).
  \end{theorem}
  \begin{proof}
    As in Theorem \ref{theo:DBI}, the proof begins by using the results from Lemma \ref{lem:DBI} and the expression $B_{j}^{-} \leq \bar{\lambda}_{j} \leq B_{j}^{+}$ 
    to bound $\bar{\lambda}_{n-2}\delta_{n-1}$, $\bar{\lambda}_{n-3}\delta_{n-2}$, $\cdots$, $\bar{\lambda}_{1}\delta_{2}$ and so to obtain the result 
    \begin{equation} \label{eq:theo_lambda_2}
      B_{1}^{-} \leq \bar{\lambda}_{1}\delta_{2} \leq B_{1}^{+}. 
    \end{equation} 
    Equation (\ref{eq:theo_lambda_2}) is then used to derive the target bounds.
    Starting with the left side of the inequality (i.e., $B_{1}^{-} \leq \bar{\lambda}_{1}\delta_{2}$)
    and noting that $\frac{1}{s(s-1)} \leq -4$ and $-\frac{1}{s} \leq -1$, yields 
    \begin{equation}
      \begin{aligned}
        \frac{m_{r}-s}{s(s-1)}d_{1}  = & \bigg(\frac{m_{r}-1}{s(s-1)}+\frac{1-s}{s(s-1)}\bigg)d_{1} \\
                                     = & \bigg(\frac{m_{r}-1}{s(s-1)}-\frac{1}{s}\bigg)d_{1} \\ 
                                  \leq & \bigg(-4(m_{r}-1)-1\bigg)d_{1} \\ 
                                     = & B_{1}^{-} \\
                                  \leq & \bar{\lambda}_{1}\delta_{2}.      
      \end{aligned}
    \end{equation}
    Isolating $m_{r}$, leads to the desired result
    \begin{equation}\label{eq:mr_bound}
        m_{r} \geq s + \frac{s(s-1)}{d_{1}}\bar{\lambda}_{1}\delta_{2} = \delta_{1} = S_{n}(x).
    \end{equation}
    Now, addressing the right side of the inequality (i.e. $B_{1}^{+} \geq \bar{\lambda}_{1}\delta_{2}$) and
    noting that $\frac{1}{s(s-1)} \leq -4$ and $-\frac{1}{s-1} \geq 1$, gives
    \begin{equation}
      \begin{aligned} 
        \frac{m_{\ell}-s}{s(s-1)}d_{1} = & \bigg(\frac{m_{\ell}}{s(s-1)} - \frac{s}{s(s-1)} \bigg)d_{1} \\
                                       = & \bigg(\frac{m_{\ell}}{s(s-1)} - \frac{1}{(s-1)} \bigg)d_{1} \\
                                    \geq & \bigg(-4 m_{\ell} + 1 \bigg)d_{1} \\
                                       = & B_{1}^{+} \\ 
                                    \geq & \lambda_{1}\delta_{2}.
      \end{aligned}
    \end{equation}
    Isolating $m_{\ell}$ leads to the desired bound
    \begin{equation}\label{eq:ml_bound}
        m_{\ell} \leq s + \frac{s(s-1)}{d_{1}}\bar{\lambda}_{1}\delta_{2} = \delta_{1} = S_{n}(x).
    \end{equation}
    The proof is concluded by combining Equations (\ref{eq:mr_bound}) and (\ref{eq:ml_bound}) to obtain
    \begin{equation}
        m_{\ell} \leq s + \frac{s(s-1)}{d_{1}}\bar{\lambda}_{1}\delta_{2} = \delta_{1} = S_{n}(x) \leq m_{r}.
    \end{equation}
  \end{proof}
  At the boundary intervals both the DBI and PPI methods construct the interpolants using a left- or right-biased stencil. 
  For the left boundary, the final stencil is built by successively appending mesh points from the right side of the of the interval $I_{i}$. 
  In the same way, the final stencil for the right boundary interval is obtained by successively appending the mesh points from the left side. 
  For the nearby boundary intervals, the stencil points selection process could reach the boundary before completing the final stencil. 
  In such a case, the remaining points are selected from the right if the left boundary is reached and from the left is the right boundary is reached. 

  \subsection{Hidden Local Extrema}
  \label{eq:edgecase}
  %%%
  The interval $I_{i}$ may contain a hidden extremum when two of three divided differences $U[x_{i-1}, x_{i}]$,  $U[x_{i+1}, x_{i}]$ and $U[x_{i+1}, x_{i+2}]$ of the neighboring intervals are of opposite signs.
  In this case, the PCHIP and DBI algorithms truncate the extremum whereas the relaxed nature of the PPI algorithm allows for a better approximation of the extremum.
  In \cite{berzins2010nonlinear}, when an extremum is detected, the ENO approach is used to construct the interpolant.
  The ENO approach may fail to recover the extremum or result in oscillations that violate the requirements for positivity and reduce the accuracy. 
  The data-bounded method in \cite{berzins2010nonlinear} is much more restrictive and does not address positivity.
  These limitations can be addressed by using a bounded positive interpolant.

  To simplify the notation, let us defined $\sigma_{i-1}$, $\sigma_{i}$ and $\sigma_{i+1}$ such that 
  \begin{equation}\label{eq:extremum}
    \sigma_{i-1} = U[x_{i-1}, x_{i}], \textrm{ }  \sigma_{i}=U[x_{i+1}, x_{i}], \textrm{ and } \sigma_{i+1}=U[x_{i+1}, x_{i+2}].
  \end{equation}
  As in \cite{berzins2010nonlinear} and \cite{sekora2009extremumpreserving}, we assume that there exists an extremum in $I_{i}$ if  
  \begin{equation}\label{eq:extremum}
   \sigma_{i-1}\sigma_{i+1} < 0, \textrm{ or } \sigma_{i-1}\sigma_{i} < 0.
  \end{equation}
  To address the cases with and without extremum, we choose the parameters $\Delta_{min}$ and $\Delta_{max}$ according to
  \begin{equation} \label{eq:umin2}
    \Delta_{min} =
    \begin{cases}
      \big|\mathbf{min}\big(u_{i}, u_{i+1}\big) \big| & \textrm{if } \sigma_{i-1} \sigma_{i+1} < 0 
                                               \textrm{ and } \sigma_{i-1} < 0 \\
                                             & \textrm{or } \sigma_{i-1}\sigma_{i+1} \geq 0 
                                               \textrm{ and } \sigma_{i-1} \sigma_{i} < 0 \\
      \epsilon \big|\mathbf{min}\big(u_{i}, u_{i+1}\big) \big| & \textrm{otherwise},
    \end{cases}
  \end{equation}
  and
  \begin{equation}\label{eq:umax2}
    \Delta_{max} = 
    \begin{cases}
      \big|\mathbf{max}\big(u_{i}, u_{i+1}\big) \big|  & \textrm{if } \sigma_{i-1} \sigma_{i+1} < 0      
                                                 \textrm{ and } \sigma_{i-1} > 0 \\             
                                              & \textrm{or } \sigma_{i-1}\sigma_{i+1} \geq 0 
                                                 \textrm{ and } \sigma_{i-1} \sigma_{i} < 0 \\ 
      \epsilon \big|\mathbf{max}\big(u_{i}, u_{i+1}\big) \big| & \textrm{otherwise}.
    \end{cases}
  \end{equation}
  $\epsilon$ is a parameter introduced to adjust $\Delta_{min}$ and $\Delta_{max}$ when no extremum is detected.
  In Equation (\ref{eq:umin2}), the interval $I_{i}$ has a local maximum if $\sigma_{i-1} \sigma_{i+1}<0$ and $\sigma_{i-1} < 0$.
  Correspondingly, in Equation (\ref{eq:umax2}), the interval $I_{i}$ has a local minimum if $\sigma_{i-1} \sigma_{i+1}<0$ and $\sigma_{i-1} > 0$.
  In both Equations (\ref{eq:umin2}) and \ref{eq:umax2}, the type of extremum is ambiguous if $\sigma_{i-1} \sigma_{i+1}$, and $\sigma_{i-1} \sigma_{i}<0$.
  When an extremum is identified, $\Delta_{min}$ and/or $\Delta_{max}$ are chosen to be sufficiently large to allow the interpolant $U^{p}(x)$ to grow beyond the data as needed to approximate the extremum without violating the requirement for positivity.
  In the case where no extremum is identified, the parameter $\epsilon$ is used to adjust $\Delta_{min}$ and/or $\Delta_{max}$ to be sufficiently large to allow higher degree interpolants compared to the DBI method, but sufficiently small to not allow for large oscillations that will degradate the accuracy of the approximation. 

\begin{figure}[H]
    \centering
    \includegraphics[scale=0.22]{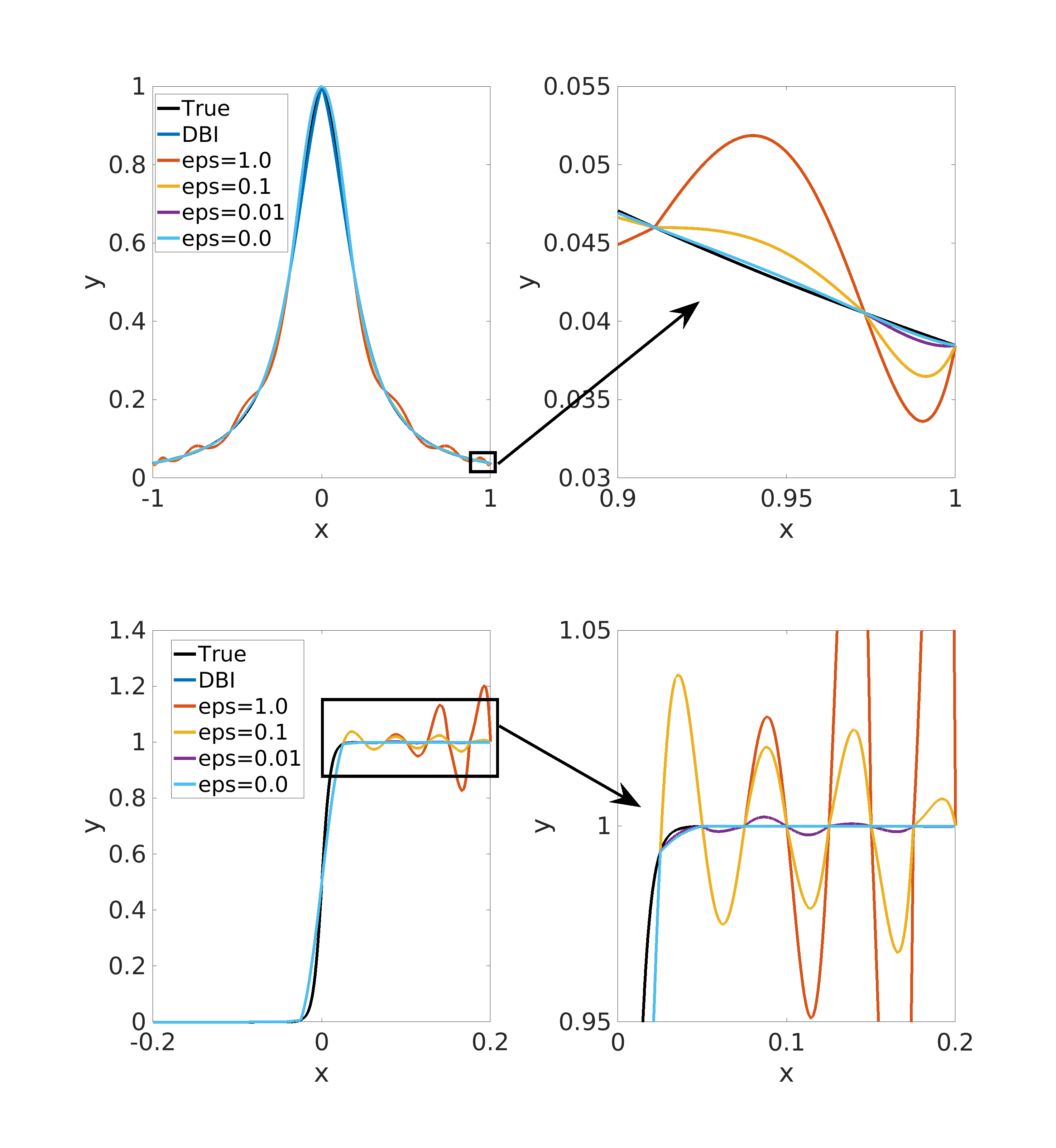}
    \caption{The top row shows an approximation of $f_{1}(x)$ from $N=17$ LGL points using DBI and PPI with different values of $\epsilon$.
      The bottom row shows an approximation of $f_{2}(x)$ from $N=17$ uniformly-spaced points using DBI and PPI with different values of $\epsilon$.
      The target polynomial degree is set to $d=16$ for both $f_{1}(x)$ and $f_{2}(x)$.}
    \label{fig:runge_oscillations}
  \end{figure}

  In Figure \ref{fig:runge_oscillations}, we approximate the Runge function with $N=17$ LGL points and different values of $\epsilon$, and the target polynomial degree is set to $d=16$ for each interval. 
  For $\epsilon > 0.01$, the PPI method leads to oscillations, whereas for $\epsilon \leq 0.01$ the oscillations are removed.
  Similar oscillations are seen when using high-order Chebyshev polynomials.
  The cutoff for the positive parameter $\epsilon$ depends on the underlying function and the input data.
  For the Runge example with $N=17$ uniformly-spaced points, the spurious oscillations are removed for $\epsilon \leq 0.05$.
  With the same Runge example with $N=129$ and $d=16$, the unconstrained approximation does not produce oscillations and $\epsilon$ can be set to any value in $[0, 1]$.
  In the case of the smoothed Heaviside examples, setting $\epsilon = 0.05$ with $N=17$ uniformly-spaced points lead to large oscillations that degrade the approximations.
  However, for $\epsilon \leq 0.01$ with $N=17$,  the oscillations are significantly reduced, and the approximation improved, as shown on the bottom part of Figure \ref{fig:runge_oscillations}.
  Setting $\epsilon =0.0$ will completely eliminate the oscillations.
  Overall, using $\epsilon \leq 0.01$ is sufficient to remove or significantly reduce the oscillations and improve the approximation. 
  For an interval $I_{i}$ with no extremum, as $\epsilon$ approaches zero and both $\Delta_{min}$ and $\Delta_{max}$ get smaller, the approximation method becomes closer to the DBI approach.
  As for the DBI approach, the PPI method may become restrictive for higher degree polynomial interpolants as $\epsilon$ approaches zero. 
  This approach is also further explored for a variety of practical applications \cite{tajo20222PPIsoftware}.

  The right part of Figure \ref{fig:runge_oscillations} shows the interpolants used at the right boundaries in both the Runge and smoothed Heaviside examples. 
  At the right boundary of the Runge example, the stencil $\{x_{N-12} \cdots x_{N}\}$ is used to build the data-bounded interpolant and the stencil $\{x_{N-16}, \cdots,  x_{N}\}$ is used for the positive interpolant with $\epsilon=1$. 
  As the positive parameter $\epsilon$ gets smaller the upper and lower bounds for the interpolant gets tighter and converges to the DBI bounds. 
  The stencil used for both the DBI and PPI are the same for $\epsilon \leq 0.01$. 
  At the boundary intervals the PPI method allows for higher degree interpolants compared to the DBI method. 
  However, these higher degree interpolants while positive may introduce oscillations that can be removed using the parameter $\epsilon$.

  For $u_{i}=u_{i+1}$, $m_{\ell}$, $m_r$ and $U_{n}(x)$ as written in Equations (\ref{eq:CaseI}), (\ref{eq:CaseII}) and (\ref{eq:UN3}) are not defined.
  The PPI algorithm addresses this limitation by re-writing $U_{n}(x)$ as 
    \begin{equation}
      U_{n}(x) = u_{i} + U[x_{1}^{l}, \cdots, x_{1}^{r}](x_{i+1}-x_{i})(x_{1}^{r}-x_{1}^{l}) S_{n}(x),
    \end{equation}
    where $S_{n}(x)$ is expressed as follows:
    \begin{equation}
      S_{n}(x) = \sum_{j=1}^{n-1}\bar{s}_{j}.
    \end{equation}
    The summation starts at $j=1$ because the linear term $\frac{u_{i+1}-u_{i}}{x_{i+1}-x_{i}}(x-x_{i}) =0$. 
    Let 
    \begin{equation}
      w = U[x_{1}^{l}, \cdots, x_{1}^{r}](x_{i+1}-x_{i})(x_{1}^{r}-x_{1}^{l}).
    \end{equation}
    $\bar{\lambda}_{j}$ in this context is defined as 
    \begin{equation}
    \bar{\lambda}_{j} = \frac{U[x_{j}^{l}, \cdots, x_{j}^{r}]}{w}
    \prod_{k=0}^{j}(x_{k}^{r}-x_{k}^{l}).
  \end{equation}

  For $u_{i}= u_{i+1}$, the parameters $m_{\ell}$ and $m_r$ are then defined according to
  \begin{enumerate}[label=\textbf{(\Roman*)}]
    \item : $U[x_{1}^{l}, \cdots, x_{1}^{r}] >0$
      \begin{equation}\label{eq:ml_mr_1}
        m_{\ell} = \mathbf{min}\bigg( 0, \frac{u_{min}-u_{i}}{w} \bigg), \textrm{ and }
        m_{r} = \mathbf{max}\bigg( 1, \frac{u_{max}-u_{i}}{w}\bigg)
      \end{equation}
    \item : $U[x_{1}^{l}, \cdots, x_{1}^{r}] <0$
      \begin{equation}\label{eq:ml_mr_2}
        m_{\ell} = \mathbf{min}\bigg( 0, \frac{u_{max}-u_{i}}{w}\bigg), \textrm{ and }
        m_{r} = \mathbf{max}\bigg( 1, \frac{u_{min}-u_{i}}{w}\bigg).
      \end{equation}

  \end{enumerate}
    For $U[x_{i}, x_{i+1}] = U[x_{1}^{l}, \cdots, x_{1}^{r}] =0$, the data $u_{i-1}$, $u_{i}$, $u_{i+1}$, and $u_{i+2}$ have the same value ($u_{i-1}=u_{i}=u_{i+1}=u_{i+2}$). 
  In this case, the algorithm approximates the function in the interval $I_{i}$ with a linear interpolant.
  For both cases $U[x_{1}^{l}, \cdots, x_{1}^{r}] <0$ and $U[x_{1}^{l}, \cdots, x_{1}^{r}] >0$,
  $B_{j}^{+}$ and $B_{j}^{-}$ remain defined as previously in Equations (\ref{eq:B_PPIplus}) and (\ref{eq:B_PPIminus}).
  Lemma \ref{lem:DBI} and Theorem \ref{theo:SN} still hold and remain unchanged.

\begin{figure}[H]
    \centering
    \includegraphics[scale=0.22]{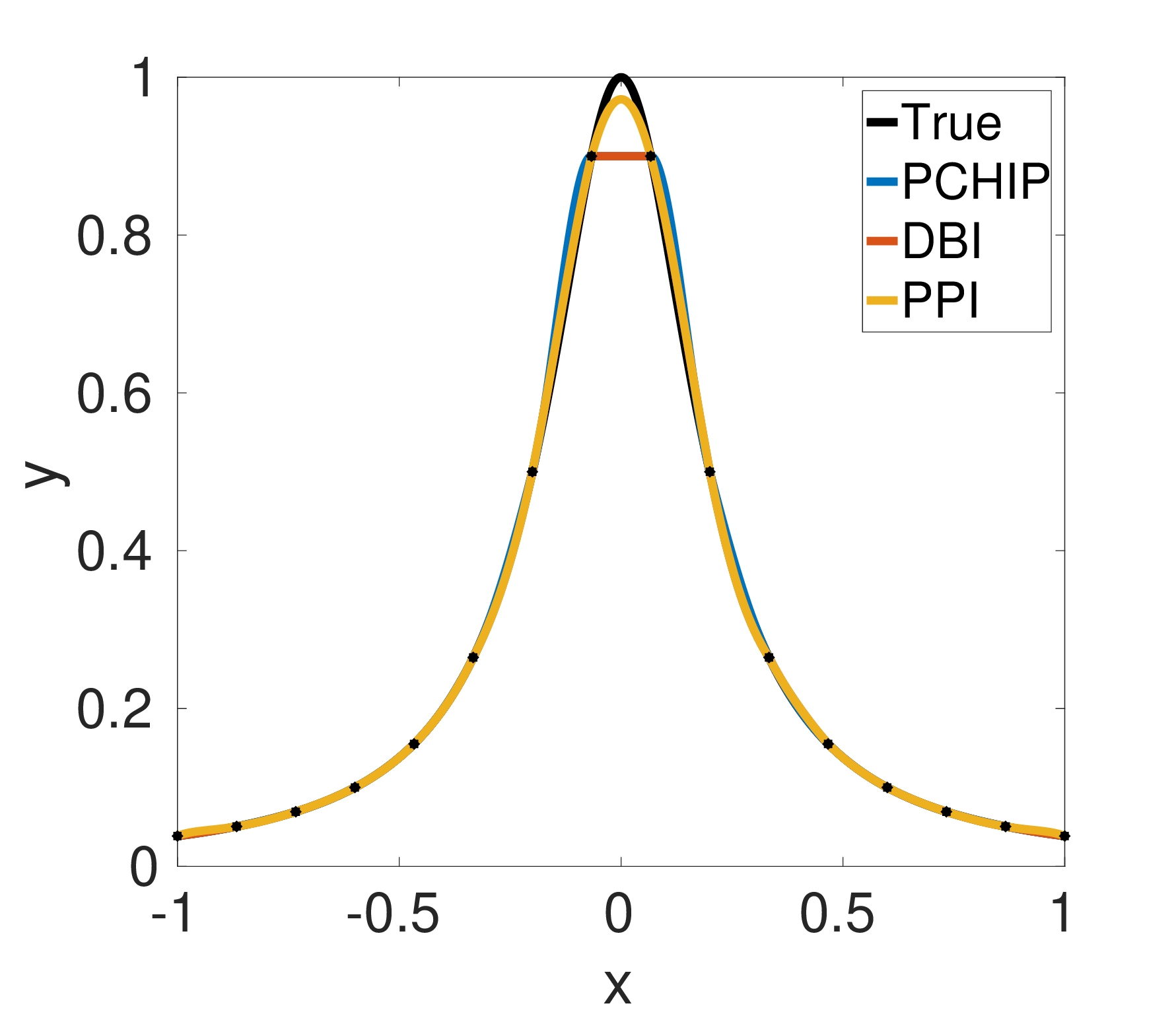}
    \caption{Approximation of $f_{1}(x)$ with $N=16$ points using PCHIP, DBI and PPI.
      The interpolants from DBI and PPI are in $\mathcal{P}_{8}$, where $8$ is the target polynomial degree.}
    \label{fig:runge}
  \end{figure}

 Figure \ref{fig:runge} shows an example with $u_{i}=u_{i+1}$ and a hidden local extremum at $x=0$. 
 In Figure \ref{fig:runge}, we approximate the Runge function $f_{1}(x)$ using the PCHIP, DBI, and PPI methods from $16$ uniformly-spaced data points.
 The PPI method is able to better capture the peak compared to the DBI and PCHIP methods.

  \subsection{Algorithm}
  \label{subsec:algo}
  The ENO reconstruction can result in a stencil that is biased to the left or right.
  Rogerson et al. \cite{rogerson1990numerical} demonstrated that a biased ENO stencil may lead to some stability issues when used to solve hyperbolic equations, and a refined resolution may lead to even larger errors.
  To address this limitation, Shu \cite{shu1990numerical} developed a modified ENO reconstruction that uses a bias coefficient to target a preferred final stencil.   
  Furthermore, a left- and right-biased stencil may fail to recover hidden local extrema.
  For instance, if $U[x_{i-1}, x_{i}]>0$, $U[x_{i}, x_{i+1}]<0$, and $U[x_{i+1}, x_{i+2}]>0$, the interval $I_{i}$ has an extremum.
  In such a case, if the points in the final stencil are all to the right or left of $x_{i}$, the interpolant may fail to recover the extremum.
  The points $x_{i-1}$ and $x_{i+2}$ are important for identifying and reconstructing a hidden local extremum.
  However, the right-biased stencils does not include $x_{i-1}$, and the left-biased stencil does not include $x_{i+2}$.
  To resolve these issues due to biased stencils, the algorithm introduced here favors a symmetric stencil over the ENO stencil in addition to enforcing the requirements for data boundedness or positivity preservation.
  A symmetric stencil centered around $x_{i}$ includes $x_{i-1}$ and $x_{i+2}$ and better approximates a hidden local extremum compared to a biased stencil.

   Before we present the algorithm for the DBI and PPI method, let us define $\bar{\lambda}_{j+1}^{-}$ and $\bar{\lambda}_{j+1}^{+}$.
   At any given step $j$, the next point inserted into $\mathcal{V}_{j}$ can be to the right or left.
   $\bar{\lambda}_{j+1}^{-}$ and $\bar{\lambda}_{j+1}^{+}$ correspond to the case where the stencil inserted is to the left and right, respectively.
   \begin{equation}
     \begin{cases}
       \bar{\lambda}_{j+1}^{-} =\bar{\lambda}_{j+1} & \textrm{ with } \mathcal{V}_{j+1} = \{x_{p} \} \cup \mathcal{V}_{j} \\ 
       \bar{\lambda}_{j+1}^{+} =\bar{\lambda}_{j+1} & \textrm{ with } \mathcal{V}_{j+1} = \mathcal{V}_{j} \cup \{x_{q} \}. 
     \end{cases}
   \end{equation}
   As a reminder, $x_{p}$ and $x_{q}$ are the mesh points immediately to the left and right of $\mathcal{V}_{j}$.
   Given $\mathcal{V}_{j}$, let $\mu_{j}^{l}$ be the number of points to the left of $x_{i}$ and $\mu_{j}^{r}$ the number of points to the right.
   Below we introduce an algorithm for DBI and PPI based on the procedures introduced above.\\
\textbf{Input:} $\{x_{i}\}_{i=0}^{n}$, $\{u_{i}\}_{i=0}^{n}$, $\{\tilde{x}_{i}\}_{i=0}^{\tilde{n}}$, $\epsilon$ and $d$.
\textbf{Output:} $\{\tilde{u}_{i}\}_{i=0}^{\tilde{n}}$.
\begin{enumerate}
    \item Select an interval $[x_{i}, x_{i+1}]$.  Let $\mathcal{V}_{0}=\{x_{i}, x_{i+1}\}=\{x_{0}^{l}, x_{0}^{r}\}$.
    \item If $\sigma_{i-1}\sigma_{i+1} < 0$ or $\sigma_{i-1}\sigma_{i} < 0$, then the interval $I_{i}$ has a hidden local extremum.
          For the boundary intervals, we assume that the divided differences to the left and right have the same sign.
    \item Compute $u_{min}$ and $u_{max}$ using Equations (\ref{eq:umin}) and (\ref{eq:umax}).
    \item Compute $m_{r}$ and $m_{\ell}$ based on Equations (50) and (51) or Equations (72) and (73). 
          For DBI, set $m_{r}=1$ and $m_{\ell} = 0$.
      \item Given a stencil $\mathcal{V}_{j}$, 
      \begin{itemize}
        \item if $B^{-}_{j+1} \leq \bar{\lambda}_{j+1}^{+} \leq B_{j+1}^{+}$ and $B^{-}_{j+1} \leq \bar{\lambda}_{j+1}^{-} \leq B_{j+1}^{+}$
        \begin{itemize}
           \item if $\mu_{j}^{l} < \mu_{j}^{r}$ then insert a new stencil point to the left;  
           \item else if $\mu_{j}^{l} > \mu_{j}^{r}$ then insert a new stencil point to the right; 
           \item else  insert a new stencil point to the right if $|\bar{\lambda}_{j+1}^{l}|\geq|\bar{\lambda}_{j+1}^{r}|$, otherwise insert a new point to left;
        \end{itemize}
        \item else if $B^{-}_{j+1} \leq \bar{\lambda}_{j+1}^{-} \leq B_{j+1}^{+}$, then insert a new stencil point to the left; 
        \item else if $B^{-}_{j+1} \leq \bar{\lambda}_{j+1}^{+} \leq B_{j+1}^{+}$, then insert a new stencil point to the right; 
      \end{itemize}
    \item This process (Steps $3$) iterates until the halting criterion that the ratio of divided differences lies outside the required bounds stated above or the stencil has $d+1$ points, with $d$ being the target degree for the interpolant. 
    \item Evaluate the final interpolant $U^{l}(x)$ (for DBI) or $U^{p}(x)$ (for PPI) at the output points $\tilde{x}_{i}$ that are in $I_{i}$.
    \item Repeat Steps $1$--$7$ for each interval in the input 1D mesh.
  \end{enumerate} 

  At the left and right boundary intervals there are no mesh points beyond the boundaries to calculate $\sigma_{i-1}$ and $\sigma_{i+1}$, respectively.
  At both boundaries $\sigma_{i-1}$ is set to $\sigma_{i+1}$ ($\sigma_{i-1}=\sigma_{i+1}$) to ensure that no new extrema are introduced. 
  At the boundary and nearby boundary intervals the algorithm allows for hidden local extrema to be recovered. 
  For example, if the right boundary interval has a hidden extremum $\sigma_{i-1} \sigma_{i} < 0$ (from Step 2) then the algorithm will relax the bounds on the interpolant and allow for the extremum to be recovered. 

\section{Numerical Experiments}
\label{sec:results}

  In this section, we present both 1D and 2D numerical experiments that demonstrate the properties of our proposed methods.
  These experimental studies use the PCHIP, DBI, and PPI methods.
  The test functions used here are taken from test problems $1$, $2$, $7$, and $10$ in \cite{arkivtajo}.
  A full suite of test problems has been undertaken by the authors in \cite{arkivtajo}. 
  In that study, nine test problems are used with both uniform and nonuniform Legendre-Gauss-Lobatto (LGL) meshes.
  The Legendre-Gauss-Lobatto mesh consists of uniform elements with eight LGL quadrature nodes \cite{hale2013fast} inside each element. 
  The number of elements is determined by $(N-1)/8$ and $(N-1)^{2}/16$ for the 1D and 2D examples.
  The integrals in the $L^{2}-$norm calculation are approximated using the trapezoid rule with $10^{4}$ and $10^{3}\times10^{3}$ uniformly-spaced points for the 1D and 2D examples, respectively. 
  The parameter $\epsilon$ is set to $0.01$ and 
  this choice is to allow the interpolant in each interval to grow beyond the data in a bounded way.
 
  For various problems, including all the examples below, a standard Lagrange interpolant leads to large oscillations and negative values. 
  While the ENO and WENO methods reduce the oscillations, they do not address the issue of preserving data boundedness or positivity.
  The DBI and PPI methods resolve both issues. 
  The numerical experiments compare the DBI and PPI methods against the widely used PCHIP method, and show approximation errors using the algorithm described in Section \ref{subsec:algo}.

  \subsection{1D Example: Runge Function}
  \label{subsec:1d_runge}
    Our first example uses the Runge \cite{doi:10.1080/00029890.1987.12000642} function, defined as follows:
    \begin{equation}\label{eq:1d_runge}
      f_{1}(x) = \frac{1}{1+25x^2}, \quad x\in [-1,1].
    \end{equation}
    Approximating the Runge function via a standard global polynomial using the set of points provided for the experiment leads to large oscillations and negative values. 

    Tables \ref{tab:f1_PCHIP} and \ref{tab:f1_DBI_PPI} show $L^{2}$-errors and convergence rates when approximating the Runge function $f_{1}(x)$ using the uniform and LGL meshes.
    For the approximations in Table \ref{tab:f1_PCHIP}, we use the PCHIP, DBI, and PPI methods with a target polynomial degree $d=3$; 
    whereas in Table \ref{tab:f1_DBI_PPI}, we use the DBI and PPI methods with the target polynomial degree varying from $d=1$ to $d=16$.
    The results in Table \ref{tab:f1_PCHIP} show that the DBI and PPI methods lead to smaller errors and larger convergence rates compared to PCHIP for $N$ larger than $17$ in both the uniform and LGL mesh examples.
    For $N=17$, the PCHIP approach leads to smaller errors.
    For higher polynomial degrees, the PPI method gives better results compared to the DBI and PCHIP, as demonstrated in Table \ref{tab:f1_DBI_PPI}.
    These results demonstrate that the PPI method is a suitable approach for interpolating data from one mesh to another when the underlying function is similar to the Runge function. 
   For $N=17$ in this example,
   the higher order terms added when going from $\mathcal{P}_{8}$ to $\mathcal{P}_{16}$ increase the $L^{2}-$error norms.
   These results indicate that resolution for $N=17$ is not sufficient to see polynomial convergence when going from $\mathcal{P}_{8}$ to $\mathcal{P}_{16}$.
   The $L^{2}-$errors norms decrease with larger values of $N$.

\begin{figure}[H]
      \centering
      \includegraphics[scale=0.22]{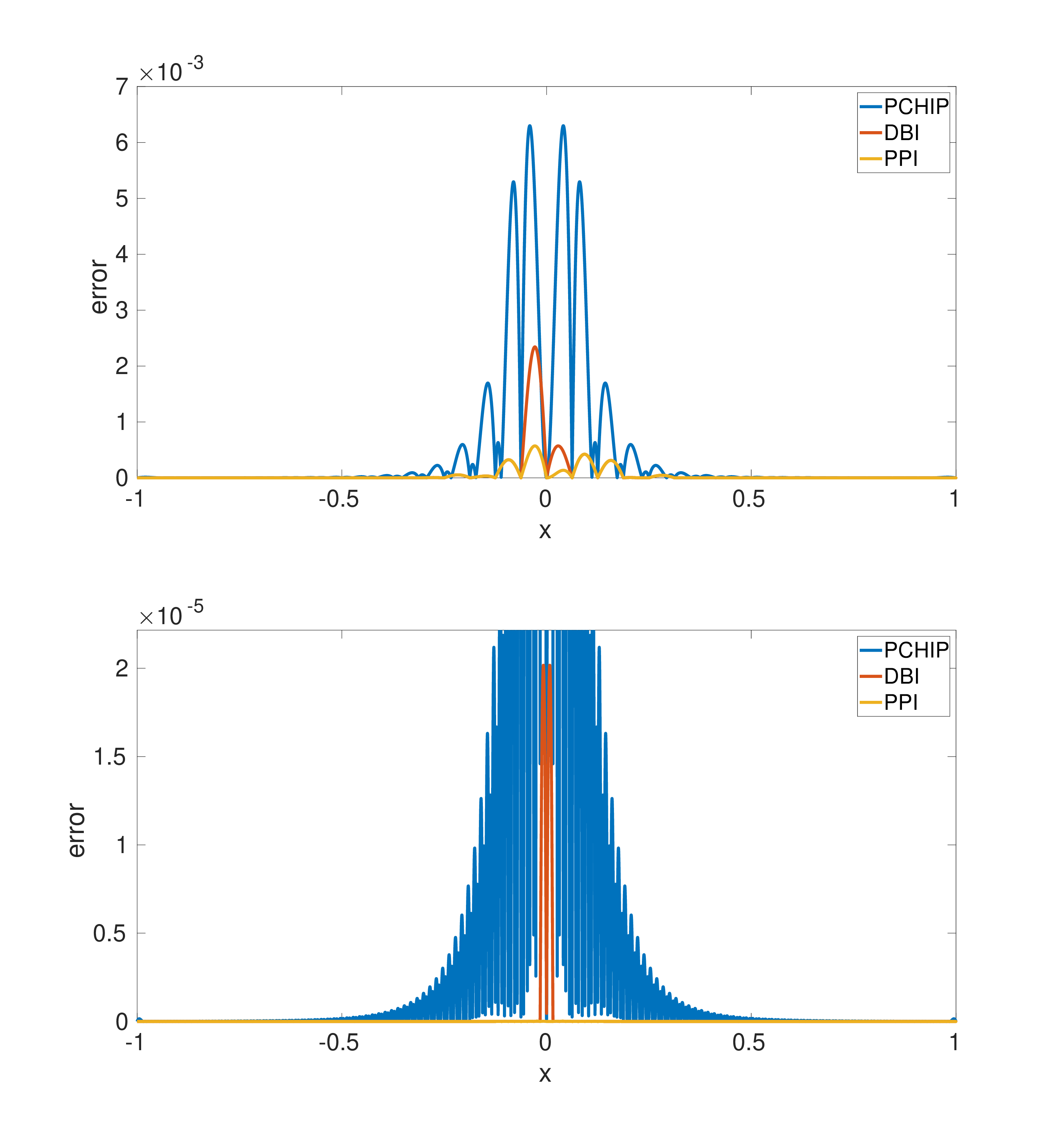}
      \caption{Error plots when approximating $f_{1}(x)$.        
      The top and bottom error plots are obtained from approximating $f_{1}(x)$ with $N=33$ and $N=129$ uniformly-spaced points, respectively. 
        The target polynomial degree is set to $d=8$ and $\epsilon =0.01$.}
      \label{fig:runge_error}
    \end{figure}

   Figure \ref{fig:runge_error} shows the errors found when approximating the Runge function $f_{1}(x)$ with PCHIP, DBI, and PPI.
   The top and bottom plots in Figure \ref{fig:runge_error} show the absolute errors when approximating the Runge example using $N=33$ and $N=129$ uniformly-spaced points, respectively.
   The target polynomial degree is set to $d=8$ for both the DBI and PPI methods and $\epsilon = 0.01$.
   The errors around the middle of the domain dominate the overall error.
   The relaxed nature of the PPI method allows for higher degree interpolants compared to the DBI and PCHIP which leads to better approximations, as shown in the bottom plots in Figure \ref{fig:runge_error}.

%    \vspace{-5.0mm}
    \begin{table}[H]
      \centering
      \begin{tabular}{ c| c c  c c c c}
        \hline
        \hline
        $N$   & PCHIP     & Rate    & DBI      & Rate     &  PPI      & Rate  \\
        \hline
                  \multicolumn{6}{c}{Uniform Mesh} \\
        \hline
         17 	 & 7.15E-03   & --     & 1.01E-02   & --     & 1.01E-02   & --     \\ 
         33 	 & 1.91E-03   & 1.99   & 1.21E-03   & 3.20   & 1.59E-03   & 2.78   \\ 
         65 	 & 3.70E-04   & 2.42   & 9.64E-05   & 3.73   & 1.12E-04   & 3.92   \\ 
         129 	 & 6.79E-05   & 2.47   & 6.29E-06   & 3.98   & 6.29E-06   & 4.20   \\ 
         257 	 & 1.22E-05   & 2.49   & 3.94E-07   & 4.02   & 3.94E-07   & 4.02   \\ 
        \hline
                  \multicolumn{6}{c}{LGL Mesh}   \\
        \hline
         17 	 & 4.75E-03   &  --    & 8.36E-03   &  --    & 8.38E-03   &  --    \\
         33 	 & 1.30E-03   & 1.96   & 1.84E-03   & 2.28   & 1.84E-03   & 2.28   \\
         65 	 & 2.86E-04   & 2.23   & 2.05E-04   & 3.24   & 2.05E-04   & 3.24   \\
         129 	 & 5.81E-05   & 2.32   & 1.17E-05   & 4.17   & 1.17E-05   & 4.17   \\
         257 	 & 1.15E-05   & 2.35   & 1.04E-06   & 3.51   & 1.04E-06   & 3.51   \\
        \hline
        \hline
      \end{tabular}                                                                                  
      \caption{$L^2$-errors and rates of convergence when using the PCHIP, DBI, and PPI  methods to approximate the function $f_{1}(x)$.
               $N$ represents the number of input points used to build the approximation.
               The approximation functions for the DBI and PPI methods are cubic interpolants.}
% in $\mathcal{P}_{3}$, where $d=3$ is the target degree.}
      \label{tab:f1_PCHIP}
    \end{table}
    \vspace{-8.0mm}
    \begin{table}[H]
      \centering
      \begin{tabular}{ c| c  c  c  c  c  c  c  c  c}
        \hline
        \hline
                  & \multicolumn{4}{c}{Uniform Mesh} && \multicolumn{4}{c}{LGL Mesh}   \\
        \hline
        $N$   & \multicolumn{2}{c}{DBI} & \multicolumn{2}{c}{PPI} && \multicolumn{2}{c}{DBI} & \multicolumn{2}{c}{PPI}   \\
        \hline
                  & $L^{2}$-error & Rate     & $L^{2}$-error & Rate   && $L^{2}$-error & Rate     & $L^{2}$-error & Rate     \\ 
        \hline
                  & \multicolumn{9}{c}{$\mathcal{P}_{1}$}               \\
        \hline
		     17 	  & 2.16E-02   &  --    & 2.16E-02   &  --    && 1.69E-02   &  --    & 1.69E-02   &  --    \\ 
		     33 	  & 6.02E-03   & 1.92   & 6.02E-03   & 1.92   && 5.84E-03   & 1.60   & 5.84E-03   & 1.60   \\ 
		     65 	  & 1.52E-03   & 2.03   & 1.52E-03   & 2.03   && 1.66E-03   & 1.86   & 1.66E-03   & 1.86   \\ 
		     129 	& 3.82E-04   & 2.02   & 3.82E-04   & 2.02   && 5.80E-04   & 1.53   & 5.80E-04   & 1.53   \\ 
		     257 	& 9.56E-05   & 2.01   & 9.56E-05   & 2.01   && 1.52E-04   & 1.94   & 1.52E-04   & 1.94   \\ 
        \hline
                  & \multicolumn{9}{c}{$\mathcal{P}_{4}$}              \\
        \hline
		     17 	  & 8.34E-03   &  --    & 7.02E-03   &  --    && 6.55E-03   &  --    & 6.54E-03   &  --    \\ 
		     33 	  & 5.91E-04   & 3.99   & 5.91E-04   & 3.73   && 7.62E-04   & 3.24   & 7.62E-04   & 3.24   \\ 
		     65 	  & 4.26E-05   & 3.88   & 2.39E-05   & 4.73   && 5.30E-05   & 3.93   & 5.29E-05   & 3.94   \\ 
		     129 	& 2.68E-06   & 4.03   & 8.00E-07   & 4.95   && 3.44E-06   & 3.99   & 3.44E-06   & 3.99   \\ 
		     257 	& 8.63E-08   & 4.99   & 2.55E-08   & 5.00   && 8.88E-08   & 5.31   & 8.87E-08   & 5.31   \\ 
        \hline
                  & \multicolumn{9}{c}{$\mathcal{P}_{8}$} \\             
        \hline
		     17 	  & 4.61E-03   &  --    & 3.11E-03   &  --    && 3.49E-03   &  --    & 4.40E-03   &  --    \\ 
		     33 	  & 4.43E-04   & 3.53   & 1.51E-04   & 4.56   && 1.76E-04   & 4.50   & 1.76E-04   & 4.85   \\ 
		     65 	  & 3.67E-05   & 3.67   & 1.05E-06   & 7.33   && 3.25E-06   & 5.89   & 3.01E-06   & 6.00   \\ 
		     129 	& 2.56E-06   & 3.88   & 3.10E-09   & 8.50   && 5.64E-08   & 5.91   & 8.82E-09   & 8.51   \\ 
		     257 	& 8.24E-08   & 4.99   & 6.80E-12   & 8.88   && 3.51E-09   & 4.03   & 3.96E-11   & 7.84   \\ 
        \hline
                  & \multicolumn{9}{c}{$\mathcal{P}_{16}$} \\
        \hline
		     17 	  & 4.34E-03   &  --    & 3.44E-03   &  --    && 4.89E-03   &  --    & 5.01E-03   &  --    \\ 
		     33 	  & 4.21E-04   & 3.52   & 4.85E-05   & 6.43   && 1.18E-04   & 5.62   & 1.17E-04   & 5.67   \\ 
		     65 	  & 3.67E-05   & 3.60   & 5.92E-08   & 9.89   && 1.22E-06   & 6.75   & 9.40E-08   & 10.51   \\ 
		     129 	& 2.56E-06   & 3.88   & 4.21E-12   & 13.94   && 5.57E-08   & 4.50   & 1.02E-11   & 13.32   \\ 
		     257 	& 8.24E-08   & 4.99   & 2.18E-16   & 14.32   && 3.51E-09   & 4.01   & 5.04E-16   & 14.38   \\ 
     \hline
     \hline
      \end{tabular}                                                                                  
      \caption{$L^2$-errors and rates of convergence when using the DBI and PPI methods to approximate the function $f_{1}(x)$.
               $N$ represents the number of input points used to build the approximation.
               The interpolants are in $\mathcal{P}_{j}$, where $j$ is the target polynomial degree.}
               %The interval $[-1,1]$ is divided into $(N_{i}-1) / j$ elements and $j+1$ quadrature points are used in each element.}
      \label{tab:f1_DBI_PPI}
    \end{table}
    %\vspace{-5.0mm}

  \subsection{1D Example: Smoothed Heaviside Function}
  \label{subsec:1d_heaviside}
    This 1D example uses an analytic approximation of the Heaviside function defined as
    \begin{equation}\label{eq:f2}
         f_{2}(x) = \frac{1}{1+ e^{-2kx}}, \quad k=100 \textrm{, and } x \in [-0.2,0.2].
    \end{equation}
    A polynomial approximation of $f_{2}(x)$ is challenging because of the large solution gradient around $x =0.$
    Attempts to use a standard polynomial approximation for this function result in oscillations and negative values.

    Tables \ref{tab:f2_PCHIP} and \ref{tab:f2_DBI_PPI} show $L^{2}$-errors and convergence rates when approximating the smoothed Heaviside function $f_{2}(x)$ using the uniform and LGL meshes.
    Table \ref{tab:f2_DBI_PPI} shows that for a target polynomial of degree $d=3$, the errors for PCHIP, DBI, and PPI are comparable.
    When the target degree increases from $d=1$ to $d=16$, the errors for the DBI and PPI methods decrease, as shown in Table \ref{tab:f2_DBI_PPI}.
    Overall, the errors from the DBI and PPI methods are comparable with DBI yielding slightly smaller errors than PPI.
    The uniform mesh leads to better approximation results compared to the LGL mesh.
    These results demonstrate that the DBI and PPI methods are both suitable for mapping data between different meshes when the underlying function is similar to the smoothed Heaviside function.

\begin{figure}[H]
      \centering
      \includegraphics[scale=0.22]{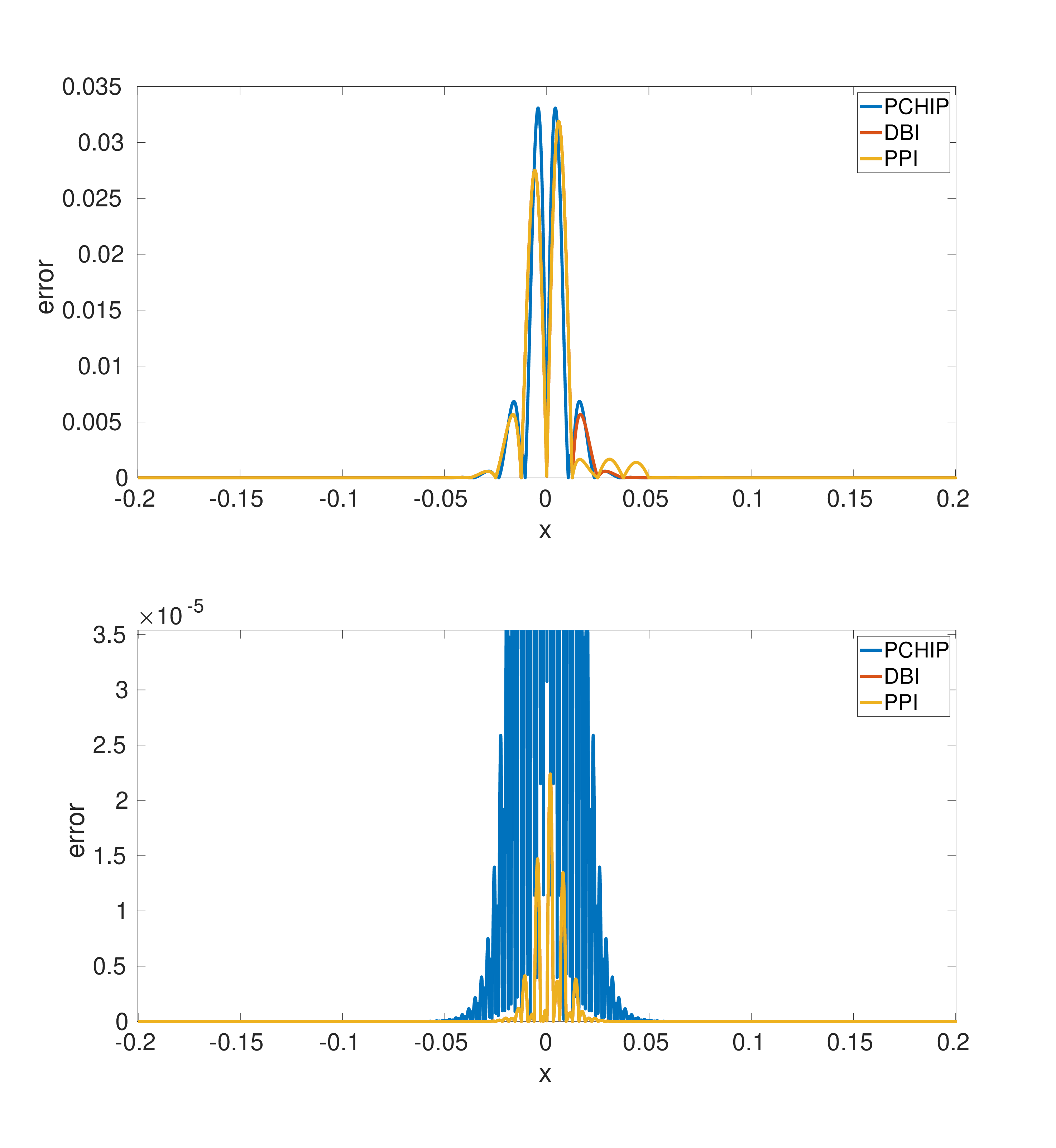}
      \caption{ Error plots when approximating $f_{2}(x)$.        
      The top and bottom error plots are obtained from approximating $f_{1}(x)$ with $N=33$ and $N=129$ uniformly-spaced points, respectively. 
        The target polynomial degree is set to $d=8$ and $\epsilon =0.01$.}
      \label{fig:heaviside_error}
    \end{figure}

    Figure \ref{fig:heaviside_error} provides examples of error plots for approximating the smoothed Heaviside function $f_{2}(x)$ with PCHIP, DBI, and PPI.
    The top and bottom plots in Figure \ref{fig:heaviside_error} show the absolute error when approximating the smoothed Heaviside function $f_{2}(x)$ using $N=33$ and $N=129$ uniformly-spaced points, respectively.
    The global error is dominated by the errors in the region with the steep gradient around $x=0$.
    The error from DBI and PPI are identical for $N=129$ because the stencil selected by both methods are the same around the region with the steep gradients.
    Away from the steep gradient the DBI and PPI methods use different stencils but the errors in those regions are negligible compared to the errors around $x=0$. 

    \begin{table}[H]
      \centering
      \begin{tabular}{ c| c  c  c  c  c c}
        \hline
        \hline
        $N$   & PCHIP     & Rate     & DBI      & Rate      &  PPI      & Rate  \\
        \hline
                  & \multicolumn{6}{c}{Uniform Mesh} \\
        \hline
	       17 	 & 2.02E-02   &  --    & 1.97E-02   &  --    & 1.97E-02   &  --   \\ 
	       33 	 & 3.38E-03   & 2.70   & 3.53E-03   & 2.59   & 3.54E-03   & 2.59  \\ 
	       65 	 & 3.59E-04   & 3.31   & 5.00E-04   & 2.88   & 5.00E-04   & 2.89  \\ 
	       129 	 & 4.21E-05   & 3.13   & 4.51E-05   & 3.51   & 4.51E-05   & 3.51  \\ 
	       257 	 & 5.12E-06   & 3.06   & 3.01E-06   & 3.93   & 3.01E-06   & 3.93  \\ 
        \hline
                  & \multicolumn{6}{c}{LGL Mesh} \\
        \hline
         17 	 & 3.65E-03   &  --    & 5.38E-03   &  --    & 5.38E-03   &  --    \\ 
         33 	 & 1.45E-03   & 1.39   & 1.55E-03   & 1.88   & 1.56E-03   & 1.86   \\ 
         65 	 & 4.07E-04   & 1.87   & 6.49E-04   & 1.28   & 6.49E-04   & 1.30   \\ 
         129 	 & 8.85E-05   & 2.23   & 9.77E-05   & 2.76   & 9.77E-05   & 2.76   \\ 
         257 	 & 1.38E-05   & 2.70   & 9.06E-06   & 3.45   & 9.06E-06   & 3.45   \\ 
        \hline
        \hline
      \end{tabular}                                                                                  
      \caption{$L^2$-errors and rates of convergence when using the PCHIP, BDI, and PPI methods to approximate the function $f_{2}(x)$.
               $N$ represents the number of input points used to build the approximation.
               The approximation functions for the DBI and PPI methods are cubic interpolants.}
               %The interval $[-0.2,0.2]$ is divided into $(N_{i}-1) / 3$ elements and $4$ quadrature points are used in each element.}
      \label{tab:f2_PCHIP}
    \end{table}
    %%%
\newpage
    \begin{table}[H]
      \centering
      \begin{tabular}{ c| c  c  c  c  c  c  c  c  c}
        \hline
        \hline
                  & \multicolumn{4}{c}{Uniform Mesh} && \multicolumn{4}{c}{LGL Mesh}   \\
        \hline
        $N$   & \multicolumn{2}{c}{DBI} & \multicolumn{2}{c}{PPI} && \multicolumn{2}{c}{DBI} & \multicolumn{2}{c}{PPI}   \\
        \hline
                  & $L^{2}$-error & Rate     & $L^{2}$-error & Rate   && $L^{2}$-error & Rate    & $L^{2}$-error & Rate     \\ 
        \hline
                  & \multicolumn{9}{c}{$\mathcal{P}_{1}$}               \\
        \hline
	       17 	  & 2.89E-02   &  --    & 2.89E-02   &  --    && 8.58E-03   &  --    & 8.58E-03   &  --    \\ 
	       33 	  & 7.69E-03   & 1.99   & 7.69E-03   & 1.99   && 5.24E-03   & 0.74   & 5.24E-03   & 0.74   \\ 
	       65 	  & 1.80E-03   & 2.14   & 1.80E-03   & 2.14   && 2.20E-03   & 1.28   & 2.20E-03   & 1.28   \\ 
	       129 	& 4.58E-04   & 2.00   & 4.58E-04   & 2.00   && 8.08E-04   & 1.47   & 8.08E-04   & 1.47   \\ 
	       257 	& 1.15E-04   & 2.00   & 1.15E-04   & 2.00   && 2.01E-04   & 2.01   & 2.01E-04   & 2.01   \\ 
        \hline
                    & \multicolumn{9}{c}{$\mathcal{P}_{4}$}              \\
        \hline
	       17 	  & 2.23E-02   &  --    & 2.23E-02   &  --    && 5.24E-03   &  --    & 5.24E-03   &  --    \\ 
	       33 	  & 4.09E-03   & 2.56   & 4.10E-03   & 2.56   && 1.10E-03   & 2.36   & 1.11E-03   & 2.34   \\ 
	       65 	  & 3.05E-04   & 3.83   & 3.05E-04   & 3.84   && 3.06E-04   & 1.88   & 3.07E-04   & 1.89   \\ 
	       129   & 1.35E-05   & 4.55   & 1.35E-05   & 4.55   && 3.32E-05   & 3.24   & 3.32E-05   & 3.24   \\ 
	       257   & 4.71E-07   & 4.87   & 4.71E-07   & 4.87   && 1.17E-06   & 4.85   & 1.17E-06   & 4.85   \\ 
        \hline
                  & \multicolumn{9}{c}{$\mathcal{P}_{8}$}     \\
        \hline
	       17 	  & 2.08E-02   &  --    & 2.08E-02   &  --    && 4.87E-03   &  --    & 4.68E-03   &  --    \\ 
	       33 	  & 3.36E-03   & 2.75   & 3.33E-03   & 2.76   && 8.71E-04   & 2.59   & 7.84E-04   & 2.69   \\ 
	       65 	  & 1.38E-04   & 4.70   & 1.38E-04   & 4.69   && 7.57E-05   & 3.60   & 1.24E-04   & 2.72   \\ 
	       129 	& 1.22E-06   & 6.90   & 1.22E-06   & 6.90   && 2.17E-06   & 5.19   & 2.17E-06   & 5.90   \\ 
	       257 	& 4.44E-09   & 8.15   & 4.44E-09   & 8.15   && 1.95E-08   & 6.83   & 1.95E-08   & 6.83   \\ 
        \hline
                 & \multicolumn{9}{c}{$\mathcal{P}_{16}$}      \\
        \hline
	       17 	  & 2.00E-02   &  --    & 2.00E-02   &  --    && 4.83E-03   &  --    & 4.64E-03   &  --    \\ 
	       33 	  & 2.93E-03   & 2.90   & 2.91E-03   & 2.91   && 7.38E-04   & 2.83   & 7.27E-04   & 2.80   \\ 
	       65 	  & 9.17E-05   & 5.11   & 9.17E-05   & 5.10   && 7.60E-05   & 3.35   & 9.41E-05   & 3.02   \\ 
	       129 	& 1.70E-07   & 9.17   & 1.70E-07   & 9.17   && 2.88E-07   & 8.14   & 2.88E-07   & 8.45   \\ 
	       257 	& 2.64E-11   & 12.73   & 2.64E-11   & 12.73   && 5.39E-11   & 12.45   & 5.39E-11   & 12.45   \\ 
        \hline
        \hline
      \end{tabular}                                                                                  
      \caption{$L^2$-errors and rates of convergence when using the DBI and PPI methods to approximate the function $f_{2}(x)$.
               $N$ represents the number of input points used to build the approximation.
               The interpolants are in $\mathcal{P}_{j}$, where $j$ is the target polynomial degree.}
      \label{tab:f2_DBI_PPI}
    \end{table}
\newpage

  \subsection{2D Example: Runge Function}
  \label{subsec:2d_runge}
    The 2D example uses an extended version of the 1D Runge function defined in Equation (\ref{eq:1d_runge}) to give
    \begin{equation}\label{eq:2d_runge}
      f_{7}(x, y) = \frac{1}{1+25(x^2+y^2)}, \quad x, y\in [-1,1].
    \end{equation}

    Tables \ref{tab:f7_PCHIP} and \ref{tab:f7_DBI_PPI} show $L^{2}$-errors and convergence rates when approximating the 2D Runge function $f_{7}(x,y)$ using the uniform and LGL meshes.
    Table \ref{tab:f7_PCHIP} compares PCHIP against DBI, and PPI with a target degree $d=3$.
    Table \ref{tab:f7_DBI_PPI} focuses on high-order interpolants using the DBI and PPI methods.
    Both the DBI and PPI methods have smaller errors compared to the PCHIP approach.
    As the target polynomial degree increases, the PPI method gives better approximation results compared to DBI and PCHIP.

    The errors in the DBI method drop more slowly than PPI when more mesh points are used because the PPI method uses higher degree interpolants compared to the DBI method. 
    The bounds on the interpolants and $\bar{\lambda}_{j}$  for data-boundedness (DBI) are more restrictive than the bounds for positivity (PPI). 
    In Tables \ref{tab:f7_DBI_PPI} and \ref{tab:f1_DBI_PPI}, when going from $d=8$ to $16$ with $N=65$, $129$, and $257$ the relaxed nature of the PPI method allows for more stencil points to be used to construct the final interpolant for each interval. 
    However, the conditions for data-boundedness are more restrictive and do not allow for more stencil points to be added when going from $d=8$ to $16$. 
    \vspace{-5mm}
    \begin{table}[H]
      \centering
      \begin{tabular}{ c| c  c  c  c  c  c }
        \hline
        \hline
        $N^{2}$   & PCHIP     & Rate    & DBI      & Rate      &  PPI      & Rate  \\
        \hline
                   \multicolumn{6}{c}{Uniform Mesh} \\
        \hline
	       $17^{2} $ & 5.01E-03   &  --    & 7.14E-03   &  --    & 7.28E-03   &  --  \\  
	       $33^{2} $ & 1.23E-03   & 2.12   & 7.82E-04   & 3.33   & 8.55E-04   & 3.23 \\  
	       $65^{2} $ & 2.33E-04   & 2.45   & 5.65E-05   & 3.88   & 5.59E-05   & 4.02 \\  
	       $129^{2}$ & 4.27E-05   & 2.48   & 3.59E-06   & 4.02   & 3.63E-06   & 3.99 \\  
	       $257^{2}$ & 7.72E-06   & 2.48   & 2.24E-07   & 4.03   & 2.27E-07   & 4.02 \\  
        \hline
                   \multicolumn{6}{c}{LGL Mesh} \\
        \hline
         $17^{2} $ & 3.26E-03   &  --    & 5.62E-03   &  --    & 5.60E-03   &  --  \\
         $33^{2} $ & 8.58E-04   & 2.01   & 1.09E-03   & 2.48   & 1.09E-03   & 2.47 \\
         $65^{2} $ & 1.88E-04   & 2.24   & 1.17E-04   & 3.29   & 1.17E-04   & 3.29 \\
         $129^{2}$ & 3.75E-05   & 2.35   & 7.05E-06   & 4.09   & 7.05E-06   & 4.09 \\
         $257^{2}$ & 7.32E-06   & 2.37   & 6.07E-07   & 3.56   & 6.07E-07   & 3.56 \\
        \hline
        \hline
      \end{tabular}                                                                                  
      \caption{$L^2$-errors and rates of convergence when using the PCHIP, DBI, and PPI methods to approximate the function $f_{7}(x, y)$.
               $N^{2}$ represents the number of input points used to build the approximation. 
               The approximation functions for the DBI and PPI methods are cubic interpolants.}
      \label{tab:f7_PCHIP}
    \end{table}
%    \vspace{-8mm}
\newpage
    \begin{table}[H]
      \centering
      \begin{tabular}{ c| c  c  c  c  c  c  c  c  c}
        \hline
        \hline
                  & \multicolumn{4}{c}{Uniform Mesh} && \multicolumn{4}{c}{LGL Mesh}   \\
        \hline
        $N^{2}$   & \multicolumn{2}{c}{DBI} & \multicolumn{2}{c}{PPI} && \multicolumn{2}{c}{DBI} & \multicolumn{2}{c}{PPI}   \\
        \hline
                  & $L^{2}$-error & Rate     & $L^{2}$-error & Rate   && $L^{2}$-error & Rate     & $L^{2}$-error & Rate     \\ 
        \hline
                  & \multicolumn{9}{c}{$\mathcal{P}_{1}$}               \\
        \hline
		     $17^{2} $ & 1.60E-02   &  --    & 1.60E-02   &  --    &&  1.10E-02   &  --    & 1.10E-02   &  --    \\
		     $33^{2} $ & 4.42E-03   & 1.94   & 4.42E-03   & 1.94   &&  3.62E-03   & 1.68   & 3.62E-03   & 1.68   \\
		     $65^{2} $ & 1.12E-03   & 2.02   & 1.12E-03   & 2.02   &&  1.20E-03   & 1.62   & 1.20E-03   & 1.62   \\
		     $129^{2}$ & 2.82E-04   & 2.02   & 2.82E-04   & 2.02   &&  4.27E-04   & 1.51   & 4.27E-04   & 1.51   \\
		     $257^{2}$ & 7.06E-05   & 2.01   & 7.06E-05   & 2.01   &&  1.11E-04   & 1.96   & 1.11E-04   & 1.96   \\
        \hline
                  & \multicolumn{9}{c}{$\mathcal{P}_{4}$}              \\
        \hline
		     $17^{2} $ & 5.07E-03   &  --    & 4.63E-03   &  --    &&  4.02E-03   &  --    & 4.05E-03   &  --    \\
		     $33^{2} $ & 3.71E-04   & 3.94   & 3.60E-04   & 3.85   &&  4.45E-04   & 3.32   & 4.45E-04   & 3.33   \\
		     $65^{2} $ & 2.62E-05   & 3.91   & 1.31E-05   & 4.89   &&  3.08E-05   & 3.94   & 3.08E-05   & 3.94   \\
		     $129^{2}$ & 1.23E-06   & 4.46   & 4.36E-07   & 4.96   &&  1.88E-06   & 4.08   & 1.88E-06   & 4.08   \\
		     $257^{2}$ & 4.96E-08   & 4.66   & 1.39E-08   & 5.00   &&  4.80E-08   & 5.32   & 4.79E-08   & 5.32   \\
        \hline                                                                   
                  & \multicolumn{9}{c}{$\mathcal{P}_{8}$}              \\       
        \hline                                                                   
		     $17^{2} $ & 3.24E-03   &  --    & 3.41E-03   &  --    &&  3.56E-03   &  --    & 3.47E-03   &  --    \\
		     $33^{2} $ & 2.88E-04   & 3.65   & 1.95E-04   & 4.31   &&  9.39E-05   & 5.48   & 9.34E-05   & 5.45   \\
		     $65^{2} $ & 2.35E-05   & 3.70   & 5.14E-07   & 8.76   &&  1.80E-06   & 5.83   & 1.51E-06   & 6.08   \\
		     $129^{2}$ & 1.16E-06   & 4.39   & 1.49E-09   & 8.53   &&  4.43E-08   & 5.41   & 4.53E-09   & 8.48   \\
		     $257^{2}$ & 4.78E-08   & 4.63   & 3.25E-12   & 8.89   &&  1.73E-09   & 4.71   & 1.87E-11   & 7.96   \\
        \hline                                                                   
                  & \multicolumn{9}{c}{$\mathcal{P}_{16}$}              \\      
        \hline                                                                   
		     $17^{2} $ & 3.69E-03   &  --    & 3.89E-03   &  --    &&  4.18E-03   &  --    & 4.18E-03   &  --    \\
		     $33^{2} $ & 2.85E-04   & 3.86   & 1.85E-04   & 4.59   &&  5.62E-05   & 6.50   & 5.68E-05   & 6.48   \\
		     $65^{2} $ & 2.35E-05   & 3.68   & 2.63E-08   & 13.07   &&  1.19E-06   & 5.69   & 4.28E-08   & 10.61   \\
		     $129^{2}$ & 1.16E-06   & 4.38   & 1.77E-12   & 14.01   &&  5.42E-08   & 4.51   & 4.31E-12   & 13.43   \\
		     $257^{2}$ & 4.76E-08   & 4.64   & 1.89E-15   & 9.93   &&  2.02E-09   & 4.77   & 1.02E-14   & 8.77   \\
        \hline
        \hline
      \end{tabular}                                                                                  
      \caption{$L^2$-errors and rates of convergence when using the DBI and PPI methods to approximate the function $f_{7}(x, y)$.
               $N^{2}$ represents the number of input points used to build the approximation.
               The interpolants are in $\mathcal{P}_{j}$, where $j$ is the target polynomial degree.}
      \label{tab:f7_DBI_PPI}
    \end{table}

  \subsection{2D Example: Smoothed Heaviside Function}
  \label{subsec:2d_heaviside}
    This 2D example uses an extension of the 1D approximation of the Heaviside function $f_{2}(x)$ defined in Equation (\ref{eq:f2}).
    The extended version is defined as follows:
    \begin{equation}\label{eq:Heaviside2}
      f_{10}(x,y) = \frac{1}{1+e^{-\sqrt{2}k(x+y)}}, \quad k=100 \textrm{, and } x, y \in [-0.2, 0.2].
    \end{equation}
    The function $f_{10}(x,y)$ is challenging because of the large gradient at $y=-x$.
 
    Tables \ref{tab:f10_PCHIP} and \ref{tab:f10_DBI_PPI} show $L^{2}$-errors and convergence rates when approximating the smoothed Heaviside function $f_{10}(x,y)$ using the uniform and LGL meshes.
    For a target polynomial of degree $d=3$, the errors for PCHIP, DBI, and PPI are comparable.
    As the target degree increases, the errors for the DBI and PPI decrease, as shown in Table \ref{tab:f2_DBI_PPI}.
    Overall, the errors from the DBI and PPI approaches are similar.

    In Tables \ref{tab:f10_DBI_PPI} and \ref{tab:f2_DBI_PPI}, the errors for the DBI and PPI methods are the same because the example used has no extrema and the interpolants used in the regions with steep gradients are the same for both the DBI and PPI methods.
    The global errors in both examples are dominated by errors in the regions with steep gradients. 
    These regions are around $x=0$ and $y=-x$ for the 1D and 2D examples, respectively.
    Away from the steep gradients DBI and PPI use different interpolants and the errors are small compared to errors around $x=0$ and $y=-x$. 
    \begin{table}[H]
      \centering
      \begin{tabular}{ c| c c  c  c  c c}
        \hline
        \hline
        $N^{2}$   & PCHIP     & Rate    & DBI      & Rate      &  PPI      & Rate  \\
        \hline
                   \multicolumn{6}{c}{Uniform Mesh} \\
        \hline
         $17^{2} $  & 8.07E-03   &  --    & 1.04E-02   &  --    & 1.05E-02   &  --  \\  
	       $33^{2} $  & 1.26E-03   & 2.80   & 2.06E-03   & 2.44   & 2.05E-03   & 2.47 \\  
	       $65^{2} $  & 1.44E-04   & 3.20   & 2.38E-04   & 3.18   & 2.38E-04   & 3.17 \\  
	       $129^{2}$  & 1.63E-05   & 3.18   & 1.64E-05   & 3.90   & 1.64E-05   & 3.90 \\  
	       $257^{2}$  & 1.94E-06   & 3.08   & 1.05E-06   & 3.99   & 1.05E-06   & 3.99 \\  
        \hline
                   \multicolumn{6}{c}{LGL Mesh} \\
        \hline
         $17^{2} $  & 1.23E-02   &  --    & 1.54E-02   &  --    & 1.56E-02   &  --  \\
         $33^{2} $  & 2.51E-03   & 2.39   & 3.86E-03   & 2.09   & 3.83E-03   & 2.11 \\
         $65^{2} $  & 3.37E-04   & 2.96   & 5.53E-04   & 2.87   & 5.53E-04   & 2.86 \\
         $129^{2}$  & 4.19E-05   & 3.04   & 4.09E-05   & 3.80   & 4.09E-05   & 3.80 \\
         $257^{2}$  & 5.96E-06   & 2.83   & 2.50E-06   & 4.05   & 2.50E-06   & 4.05 \\
        \hline
        \hline
      \end{tabular}                                                                                  
      \caption{$L^2$-errors and rates of convergence when using the PCHIP, DBI, and PPI methods to approximate the function $f_{10}(x, y)$.
               $N^{2}$ represents the number of input points used to build the approximation.
               The approximation functions for the DBI and PPI methods are cubic interpolants.}

      \label{tab:f10_PCHIP}
    \end{table}
    %% 
%%    \vspace{-10mm}
\newpage
    \begin{table}[H]
      \centering
      \begin{tabular}{ c| c  c  c  c  c  c  c  c  c  c}
        \hline
        \hline
                  & \multicolumn{4}{c}{Uniform Mesh} && \multicolumn{4}{c}{LGL Mesh}   \\
        \hline
        $N^{2}$   & \multicolumn{2}{c}{DBI} & \multicolumn{2}{c}{PPI} && \multicolumn{2}{c}{DBI} & \multicolumn{2}{c}{PPI}   \\
        \hline
                  & $L^{2}$-error & Rate     & $L^{2}$-error & Rate   && $L^{2}$-error & Rate     & $L^{2}$-error & Rate     \\ 
        \hline
                  & \multicolumn{9}{c}{$\mathcal{P}_{1}$}               \\
        \hline
		     $17^{2} $ & 1.50E-02   &  --    & 1.50E-02   &  --    &&  2.05E-02   &  --    & 2.05E-02   &  --    \\
		     $33^{2} $ & 4.57E-03   & 1.79   & 4.57E-03   & 1.79   &&  6.79E-03   & 1.66   & 6.79E-03   & 1.66   \\
		     $65^{2} $ & 1.26E-03   & 1.90   & 1.26E-03   & 1.90   &&  1.89E-03   & 1.89   & 1.89E-03   & 1.89   \\
		     $129^{2}$ & 3.23E-04   & 1.98   & 3.23E-04   & 1.98   &&  4.86E-04   & 1.98   & 4.86E-04   & 1.98   \\
		     $257^{2}$ & 8.15E-05   & 2.00   & 8.15E-05   & 2.00   &&  1.24E-04   & 1.98   & 1.24E-04   & 1.98   \\
        \hline
                  & \multicolumn{9}{c}{$\mathcal{P}_{4}$}              \\
        \hline
		     $17^{2} $ & 9.45E-03   &  --    & 9.42E-03   &  --    &&  1.37E-02   &  --    & 1.36E-02   &  --    \\
		     $33^{2} $ & 1.33E-03   & 2.95   & 1.31E-03   & 2.98   &&  2.72E-03   & 2.43   & 2.71E-03   & 2.44   \\
		     $65^{2} $ & 9.29E-05   & 3.93   & 9.29E-05   & 3.90   &&  2.39E-04   & 3.59   & 2.39E-04   & 3.58   \\
		     $129^{2}$ & 3.67E-06   & 4.71   & 3.67E-06   & 4.71   &&  1.10E-05   & 4.49   & 1.10E-05   & 4.49   \\
		     $257^{2}$ & 1.21E-07   & 4.95   & 1.21E-07   & 4.95   &&  3.90E-07   & 4.84   & 3.90E-07   & 4.84   \\
        \hline
                  & \multicolumn{9}{c}{$\mathcal{P}_{8}$}              \\
        \hline
		     $17^{2} $ & 8.04E-03   &  --    & 8.00E-03   &  --    &&  1.22E-02   &  --    & 1.21E-02   &  --    \\
		     $33^{2} $ & 1.03E-03   & 3.10   & 9.30E-04   & 3.25   &&  1.76E-03   & 2.91   & 1.75E-03   & 2.92   \\
		     $65^{2} $ & 4.83E-05   & 4.51   & 4.89E-05   & 4.35   &&  4.98E-05   & 5.26   & 4.98E-05   & 5.25   \\
		     $129^{2}$ & 2.57E-07   & 7.64   & 2.57E-07   & 7.66   &&  4.03E-07   & 7.03   & 4.03E-07   & 7.03   \\
		     $257^{2}$ & 5.27E-10   & 8.98   & 5.27E-10   & 8.98   &&  1.21E-09   & 8.42   & 1.21E-09   & 8.42   \\
        \hline
                  & \multicolumn{9}{c}{$\mathcal{P}_{16}$}              \\
        \hline
		     $17^{2} $ & 7.32E-03   &  --    & 7.31E-03   &  --    &&  1.17E-02   &  --    & 1.16E-02   &  --    \\
		     $33^{2} $ & 1.03E-03   & 2.96   & 8.90E-04   & 3.17   &&  1.46E-03   & 3.14   & 1.44E-03   & 3.15   \\
		     $65^{2} $ & 2.13E-04   & 2.32   & 2.09E-04   & 2.13   &&  1.83E-04   & 3.06   & 1.64E-04   & 3.20   \\
		     $129^{2}$ & 1.03E-06   & 7.78   & 1.03E-06   & 7.76   &&  2.15E-07   & 9.85   & 2.15E-07   & 9.69   \\
		     $257^{2}$ & 4.41E-11   & 14.59   & 4.41E-11   & 14.59   &&  9.37E-12   & 14.57   & 9.37E-12   & 14.57   \\
        \hline
        \hline
  \end{tabular}                                                                                  
      \caption{$L^2$-errors and rates of convergence when using the DBI and PPI methods to approximate the function $f_{10}(x, y)$.
               $N^{2}$ represents the number of input points used to build the approximation.
               The interpolants are in $\mathcal{P}_{j}$, where $j$ is the target polynomial degree.}
      \label{tab:f10_DBI_PPI}
    \end{table}

    \subsection{Hidden Local Extrema Examples}
    \label{subsec:extremum}
    This numerical study demonstrates the ability of the PPI method to recover hidden extrema.
    The study uses the Runge functions $f_{1}(x)$ and $f_{7}(x,y)$ with \st{a} uniform meshes.   
    The uniformly-spaced mesh points are constructed such that the extremum at $x=0$ lies inside of an interval.
    Tables \ref{tab:f1_DBI_PPI_peak} and \ref{tab:f7_DBI_PPI_peak} show $L^{2}$-error norms and convergence rates when approximating $f_{1}(x)$ and $f_{7}(x,y)$ from Equations (\ref{eq:1d_runge}) and (\ref{eq:2d_runge}).    The results from both tables show that the PPI method leads to smaller errors and larger convergence rates compared to the DBI method.
    The DBI approach uses a bounded interpolant that fails to represent the extremum at $x=0$, whereas the relaxed nature of the PPI approach allows for a more accurate representation of the extremum.
    In the case of DBI, as the target polynomial degree increases from $\mathcal{P}_{4}$ to $\mathcal{P}_{16}$, the errors and convergence rates do not improve because the global error is dominated by the local error in the interval with the hidden extremum.
    The  DBI approach only achieves an $O(h^{2.5})$ accuracy as opposed to the PPI method, that achieves the same high accuracy regardless of whether or not the extremal values are data points.
    These results highlight the advantage of the PPI method over the DBI method for recovering hidden extrema from data.
    Overall, the PPI method achieves high-order accuracy when approximating the Runge functions from data with and without hidden extrema.
    \begin{table}[H]
      \centering
      \begin{tabular}{ c| c  c   c   c  c }
        \hline
        \hline
        $N$   & \multicolumn{2}{c}{DBI} && \multicolumn{2}{c}{PPI}  \\
        \hline
                  & $L^{2}$-error & Rate     && $L^{2}$-error & Rate        \\ 
        \hline
                  & \multicolumn{5}{c}{$\mathcal{P}_{1}$}       \\
        \hline
	       16   & 2.81E-02   &  --    && 2.81E-02   &  --  \\  
	       32   & 6.41E-03   & 2.13   && 6.41E-03   & 2.13 \\  
	       64   & 1.57E-03   & 2.03   && 1.57E-03   & 2.03 \\  
	       128  & 3.88E-04   & 2.02   && 3.88E-04   & 2.02 \\  
	       256  & 9.63E-05   & 2.01   && 9.63E-05   & 2.01 \\  
        \hline
                  & \multicolumn{5}{c}{$\mathcal{P}_{4}$}      \\
        \hline
	       16   & 2.81E-02   &  --    && 1.37E-02   &  --  \\  
	       32   & 4.72E-03   & 2.57   && 6.85E-04   & 4.32 \\  
	       64   & 8.14E-04   & 2.54   && 2.57E-05   & 4.73 \\  
	       128  & 1.42E-04   & 2.52   && 8.32E-07   & 4.95 \\  
	       256  & 2.49E-05   & 2.51   && 2.60E-08   & 5.00 \\  
        \hline
                  & \multicolumn{5}{c}{$\mathcal{P}_{8}$}      \\
        \hline
	       16 	 & 2.74E-02   &  --    && 1.07E-02   &  --  \\  
	       32 	 & 4.69E-03   & 2.55   && 2.06E-04   & 5.70 \\  
	       64 	 & 8.14E-04   & 2.53   && 1.19E-06   & 7.43 \\  
	       128	 & 1.42E-04   & 2.52   && 3.32E-09   & 8.49 \\  
	       256	 & 2.49E-05   & 2.51   && 7.04E-12   & 8.88 \\  
        \hline
                  & \multicolumn{5}{c}{$\mathcal{P}_{16}$}     \\    
        \hline
	       16   & 2.75E-02   &  --    && 1.02E-02   &  --  \\  
	       32   & 4.69E-03   & 2.55   && 1.43E-04   & 6.16 \\  
	       64   & 8.14E-04   & 2.53   && 7.18E-08   & 10.96 \\  
	       128  & 1.42E-04   & 2.52   && 4.74E-12   & 13.89 \\  
	       256  & 2.49E-05   & 2.51   && 2.77E-16   & 14.06 \\  
        \hline
        \hline
      \end{tabular}                                                                              
      \caption{$L^2$-errors and rates of convergence when using the DBI and PPI methods to approximate the function $f_{1}(x)$.
               The uniform mesh used to build the approximation is constructed with $N$ points .
               The interpolants are in $\mathcal{P}_{j}$, where $j$ is the target polynomial degree.}
      \label{tab:f1_DBI_PPI_peak}
    \end{table}
\newpage
    \begin{table}[H]
      \centering
      \begin{tabular}{ c| c  c   c   c  c }
        \hline
        \hline
        $N^{2}$   & \multicolumn{2}{c}{DBI} && \multicolumn{2}{c}{PPI}  \\
        \hline
                  & $L^{2}$-error & Rate     && $L^{2}$-error & Rate        \\ 
        \hline
                  & \multicolumn{5}{c}{$\mathcal{P}_{1}$}       \\
        \hline
		     $16^{2} $ & 1.97E-02   &  --    && 1.97E-02   &  --    \\
		     $32^{2} $ & 4.71E-03   & 2.07   && 4.71E-03   & 2.07   \\
		     $64^{2} $ & 1.16E-03   & 2.02   && 1.16E-03   & 2.02   \\
		     $128^{2}$ & 2.86E-04   & 2.02   && 2.86E-04   & 2.02   \\
		     $256^{2}$ & 7.11E-05   & 2.01   && 7.11E-05   & 2.01   \\
        \hline
                  & \multicolumn{5}{c}{$\mathcal{P}_{4}$}      \\
        \hline
		     $16^{2} $ & 1.91E-02   &  --    && 7.97E-03   &  --    \\
		     $32^{2} $ & 3.27E-03   & 2.55   && 3.83E-04   & 4.38   \\
		     $64^{2} $ & 5.43E-04   & 2.59   && 1.41E-05   & 4.76   \\
		     $128^{2}$ & 9.18E-05   & 2.56   && 4.53E-07   & 4.96   \\
		     $256^{2}$ & 1.59E-05   & 2.53   && 1.41E-08   & 5.00   \\
        \hline
                  & \multicolumn{5}{c}{$\mathcal{P}_{8}$}      \\
        \hline
		     $16^{2} $ & 1.90E-02   &  --    && 6.03E-03   &  --    \\
		     $32^{2} $ & 3.27E-03   & 2.53   && 1.05E-04   & 5.84   \\
		     $64^{2} $ & 5.43E-04   & 2.59   && 5.83E-07   & 7.49   \\
		     $128^{2}$ & 9.18E-05   & 2.56   && 1.59E-09   & 8.51   \\
		     $256^{2}$ & 1.59E-05   & 2.53   && 3.36E-12   & 8.89   \\
        \hline
                  & \multicolumn{5}{c}{$\mathcal{P}_{16}$}     \\    
        \hline
		     $16^{2} $ & 1.91E-02   &  --    && 6.06E-03   &  --    \\
		     $32^{2} $ & 3.28E-03   & 2.54   && 8.11E-05   & 6.22   \\
		     $64^{2} $ & 5.43E-04   & 2.60   && 3.19E-08   & 11.31   \\
		     $128^{2}$ & 9.18E-05   & 2.56   && 2.00E-12   & 13.96   \\
		     $256^{2}$ & 1.59E-05   & 2.53   && 2.83E-15   & 9.46   \\ 
       \hline
       \hline
      \end{tabular}                                                                                  
      \caption{$L^2$-errors and rates of convergence when using the DBI and PPI methods to approximate the function $f_{7}(x, y)$.
               The uniform mesh used to build the approximation is constructed with $N^{2}$ points.
               The interpolants are in $\mathcal{P}_{j}$, where $j$ is the target polynomial degree.}
      \label{tab:f7_DBI_PPI_peak}
    \end{table}

\section{Summary and Conclusions}
\label{sec:conclusions}

In this paper, we present both an algorithm and theoretical foundations for 
sufficient conditions to ensure data boundedness and positivity on any set of mesh points via a Newton polynomial formulation.
The one-dimensional PPI and DBI methods analyzed herein are building blocks that have been extended to multidimensional PPI and DBI methods using tensor-products.
This extension consists of successively applying the one-dimensional PPI or DBI method on each dimension to generate the multidimensional results. 

The DBI method imposes restrictions on the ratio of divided differences to ensure that the interpolants are bounded by the input data.
The proof of the DBI approach presents new challenges because the configuration of mesh points may not exhibit a regular structure.
The PPI method starts from the DBI method and relaxes the bounds on the ratio of divided differences, thereby allowing the interpolants to grow beyond the data as needed while remaining positive.
The positive interpolant is further bounded by the parameters $u_{min}$ and $u_{max}$ to remove undesirable oscillations that may potentially degrade the approximation.
The proofs of both the DBI and PPI approaches rely on the results from Lemma \ref{lem:DBI}, which consist of using the definition of $B^{+}_{j}$, $B_{j}^{-}$ to arrive at the bounds $B_{j-1}^{-} \leq \bar{\lambda}_{j-1}\delta_{j} \leq B_{j-1}^{+}$.
The proofs from Theorems \ref{theo:DBI} and \ref{theo:SN} use Lemma \ref{lem:DBI} to show that $0 \leq S_{n}(x) \leq 1$ for the DBI method and $m_{\ell} \leq S_{n}(x) \leq m_r$ for the PPI method. 
 
Note that one observation we have made is that the PPI method uses higher order interpolants compared to the DBI method. 
Relaxing the bounds on the ratio of divided differences increases the range of polynomial degrees that meet the desired requirement.
The 1D and 2D numerical results, in Tables \ref{tab:f1_PCHIP}-\ref{tab:f10_DBI_PPI}, indicate that the DBI or PPI methods provided herein are appropriate for ensuring data boundedness or positivity preservation, and both methods converge as the interpolant degree and resolution increase.
Figure \ref{fig:runge_oscillations} demonstrates that enforcing positivity alone may not be sufficient to remove large oscillations.
We resolve this issue by bounding the positive polynomial with $u_{min}$ and $u_{max}$, which are determined based on user-supplied values, such as $\epsilon = 0.01$ for the numerical examples in Section \ref{sec:results}. 
In addition, Figure \ref{fig:runge} demonstrates that for an interval $I_{i}$ where there exists a local extremum, the PCHIP and DBI methods truncate the extremum whereas the PPI method leads to a better approximation of the extremum. 
The different results demonstrated that the PPI method is able to produce high-order accurate approximations in examples with and without a hidden extremum.

As this work continues, we plan to investigate different methods for accelerating the algorithm. 
The performance optimization will focus on different strategies to enable data locality and vectorization of the PPI and DBI algorithm to better take advantage of different computational architectures.
In addition, we will evaluate the use of both DBI and PPI methods for various practical applications.
This work is ongoing ~\cite{tajo20222PPIsoftware}.

\begin{acknowledgements}
This work has been supported by the US Naval Research Laboratory (559000669), 
the National Science Foundation (1521748), and the Intel Graphics and Visualization Institute at the University of Utah's Scientific Computing and Imaging (SCI) Institute (29715).
The authors would like like to thank  Dr. Alex Reinecke of the Naval Research Laboratory for his constant support and help.
\end{acknowledgements}


\bibliographystyle{spmpsci}      % mathematics and physical sciences
\bibliography{references}        % name your BibTeX data base

\end{document}